\documentclass[12pt,leqno,a4paper]{amsart}
\usepackage{amssymb,enumerate}
\newcommand{\FF}{{\mathbb{F}}}
\newcommand{\QQ}{{\mathbb{Q}}}

\newcommand{\bC}{{\mathbf{C}}}
\newcommand{\bG}{{\mathbf{G}}}
\newcommand{\bH}{{\mathbf{H}}}
\newcommand{\bL}{{\mathbf{L}}}
\newcommand{\bN}{{\mathbf{N}}}
\newcommand{\bS}{{\mathbf{S}}}
\newcommand{\bT}{{\mathbf{T}}}
\newcommand{\bU}{{\mathbf{U}}}
\newcommand{\bZ}{{\mathbf{Z}}}

\newcommand{\cC}{{\mathcal{C}}}
\newcommand{\cE}{{\mathcal{E}}}
\newcommand{\cS}{{\mathcal{S}}}

\newcommand{\fA}{{\mathfrak{A}}}
\newcommand{\fS}{{\mathfrak{S}}}

\newcommand{\Aut}{{\operatorname{Aut}}}
\newcommand{\Irr}{{\operatorname{Irr}}}
\newcommand{\GAP}{{\sf{GAP}}}

\newcommand{\GL}{{\operatorname{GL}}}
\newcommand{\PGL}{{\operatorname{PGL}}}
\newcommand{\SL}{{\operatorname{SL}}}
\newcommand{\PSL}{{\operatorname{L}}}  
\newcommand{\GU}{{\operatorname{GU}}}
\newcommand{\SU}{{\operatorname{SU}}}
\newcommand{\PSU}{{\operatorname{U}}}  
\newcommand{\Sp}{{\operatorname{Sp}}}
\newcommand{\PSp}{{\operatorname{S}}}  
\newcommand{\GO}{{\operatorname{GO}}}
\newcommand{\SO}{{\operatorname{SO}}}
\newcommand{\PSO}{{\operatorname{O}}}  
\newcommand{\Chevie}{{\sf Chevie}{}}

\newcommand{\tw}[1]{{}^{#1}\!}
\newcommand{\hbS}{\hat\bS}
\newcommand\subn{{\,\triangleleft\triangleleft\,}}
\newcommand{\Sub}{{\operatorname{Sub}}}

\let\al=\alpha
\let\eps=\epsilon

\let\Om=\Omega

\theoremstyle{theorem}
\newtheorem{thm}{Theorem}[section]
\newtheorem{lem}[thm]{Lemma}
\newtheorem{prop}[thm]{Proposition}
\newtheorem{cor}[thm]{Corollary}

\newtheorem{thmA}{Theorem}
\newtheorem{conjA}[thmA]{Conjecture}

\theoremstyle{definition}
\newtheorem{rem}[thm]{Remark}

\begin{document}

\title[Subnormalisers of semisimple elements]{Subnormalisers of semisimple elements\\ in finite groups of Lie type}

\author{Gunter Malle}
\address{FB Mathematik, RPTU Kaiserslautern, Postfach 3049,
  67653 Kaisers\-lautern, Germany.}
\email{malle@mathematik.uni-kl.de}

\begin{abstract}
We determine subnormalisers of semisimple elements of prime power order in
finite quasi-simple groups of Lie type. For this, we determine the maximal
overgroups of normalisers of Sylow tori. This is motivated by the recent
character correspondence conjecture by Moret\'o and Rizo as well as by the
question of existence of quasi-semiregular elements in finite permutation
groups.
\end{abstract}



\thanks{The author gratefully acknowledges financial support by the DFG---
 Project-ID 286237555 -- TRR 195. He also thanks Mandi Schaeffer Fry and
 Juan Mart\'inez Madrid for their comments on an earlier version.}

\date{\today}

\maketitle


\section{Introduction}

In this paper we continue our investigation in \cite{Ma25} in relation to the
recent conjecture of Moret\'o and Rizo \cite{MR25} on character correspondences
for finite groups $G$. For a prime~$p$, these should relate the irreducible
characters of $G$ to those of subnormalisers of $p$-elements of~$G$, which are
defined as follows. For a subgroup $H$ of $G$ let
$$S_G(H):=\{g\in G\mid H\subn\langle g,H\rangle \},$$
the set of elements $g\in G$ such that $H$ is subnormal in $\langle g,H\rangle$,
and define
$$\Sub_G(H):=\big\langle S_G(H)\big\rangle$$
to be the \emph{subnormaliser} of $H$. If $H=\langle x\rangle$ is generated by
a single element, we also write $\Sub_G(x)$ for $\Sub_G(\langle x\rangle)$.
For an element $x\in G$ let $\Irr^x(G)$ denote the set of complex irreducible
characters of~$G$ that do not vanish at~$x$. The following was put forward in
\cite{MR25}:

\begin{conjA}[Moret\'o--Rizo]   \label{conj:AN}
 Let $G$ be a finite group and $p$ a prime. Then for any $p$-element $x\in G$
 there exists a bijection $f_x:\Irr^x(G)\to\Irr^x(\Sub_G(x))$ such that
 \begin{enumerate}
  \item[\rm(1)] $\chi(1)_p=f_x(\chi)(1)_p$, and
  \item[\rm(2)] $\QQ(\chi(x))=\QQ(f_x(\chi)(x))$.
 \end{enumerate}
\end{conjA}

In order to investigate this conjecture for non-abelian simple groups it seems
useful to understand the structure of subnormalisers of $p$-elements. In our
predecessor paper \cite{Ma25} we classified semisimple \emph{picky}
$p$-elements, that is, elements whose subnormaliser is a Sylow $p$-normaliser,
of quasi-simple groups of Lie type except for $p\le3$, and obtained partial
information on subnormalisers of unipotent elements. The picky semisimple 2-
and 3-elements were then determined in \cite{MS25}. Here, we continue the
investigation of subnormalisers for semisimple $p$-elements in groups of Lie
type. This naturally leads to the question of understanding overgroups of
certain torus normalisers which might be of independent interest.
\smallskip

Our main result is the determination of maximal overgroups of normalisers of
Sylow $d$-tori, for integers $d\ge1$, which in turn yields the maximal
overgroups of normalisers of abelian Sylow subgroup, and using this the
description of subnormalisers of semisimple elements
lying in these abelian Sylow subgroups. It turns out, \emph{a posteriori}, that
these subnormalisers are natural geometrically defined subgroups except for
one single case in the exceptional group $G_2(4)$ where the sporadic simple
group $J_2$ of Janko appears. Combining our analysis with results for symmetric
and sporadic groups we can state:

\begin{thmA}   \label{thm:main}
 Let $G$ be a finite quasi-simple group and $p$ a prime such that $G$ has
 abelian Sylow $p$-subgroups. Then the subnormalisers of all $p$-elements
 in $G$ are known.
\end{thmA}

The situation for non-abelian Sylow subgroups is considerably more difficult
with further types of subnormalisers appearing, and not so tightly related to
Sylow tori, so we will not discuss it here. 
\smallskip

Our investigations are related to other current research work. Subnormalisers
play a central role in the investigation of Giudici, Morgan and Praeger on
finite permutation groups $G$ containing quasi-semiregular elements, since
$x\in G$ of prime order is quasi-semiregular if and only if there is a
point-stabiliser of $G$ containing $\Sub_G(x)$ (see \cite[Thm~3.3]{GMP}).

Baumeister, Burness, Guralnick and Tong-Viet \cite{BBGT} classify finite
non-abelian almost simple groups with a Sylow $p$-subgroup contained in a
unique maximal subgroup. This is related to our work for simple groups of Lie
type $G$ with an abelian Sylow $p$-subgroup~$P$ as follows. If $P$ is contained
in a unique maximal subgroup~$M$ of~$G$, then so is its normaliser $\bN_G(P)$;
in particular if this normaliser lies in several maximal overgroups, the same
is true \emph{a fortiori} for $P$. That is, all examples in \cite{BBGT} with
abelian Sylow subgroups also show up as part of our classification.
\medskip

Our paper is built up as follows.
In Section~\ref{sec:ab} we collect some basic facts on subnormalisers of
semisimple elements lying in abelian Sylow subgroups of finite groups of Lie
type.
In Section~\ref{sec:exc} we determine in Theorem~\ref{thm:over} the maximal
overgroups of normalisers of Sylow $d$-tori in exceptional groups of simply
connected Lie type and use this to describe in Theorem~\ref{thm:subn exc} the
subnormalisers of semisimple $p$-elements in these groups in the case that
Sylow $p$-subgroups are abelian. In Section~\ref{sec:class} we solve the
same problems for the various series of classical groups of Lie type, see in
particular Theorems~\ref{thm:subn SL}, \ref{thm:subn SU}, \ref{thm:subn Sp},
\ref{thm:subn SOodd}, \ref{thm:subn SOeven+} and~\ref{thm:subn SOeven-}. In
Sections~\ref{sec:symm} and~\ref{sec:spor} we complement our results by
describing the subnormalisers of $p$-elements in symmetric as well as in
sporadic simple groups with abelian Sylow $p$-subgroups and thus complete the
proof of Theorem~\ref{thm:main}.

\section{Subnormalisers of semisimple elements in the abelian Sylow case}   \label{sec:ab}

Let $\bG$ be a simple algebraic group of simply connected type with a
Frobenius endomorphism $F$ with respect to an $\FF_q$-structure and let
$\ell$ be a prime. When Sylow $\ell$-sub\-groups of $G:=\bG^F$ are abelian the
determination of subnormalisers relies on the knowledge of possible overgroups
of normalisers of Sylow $d$-tori of $\bG$. See \cite[\S24]{MT11} or
\cite[\S3.5]{GM20} for background on Sylow tori.

\begin{prop}   \label{prop:Burnside}
 Let $\ell$ be a prime not dividing $q$ such that Sylow $\ell$-subgroups of
 $G=\bG^F$ are abelian. Let $d:=e_\ell(q)$ be the order of $q$ modulo $\ell$
 and $\bS_d$ be a Sylow $d$-torus of $\bG$. Then we have
 \begin{enumerate}[\rm(a)]
  \item $\ell>2$, $\ell$ is good for $\bG$ and not a torsion prime;
  \item $\Phi_d$ is the unique cyclotomic polynomial dividing the generic order
   of $(\bG,F)$ with $\ell|\Phi_d(q)$;
  \item $\bS_d^F$ contains a Sylow $\ell$-subgroup $P$ of $G$; and
  \item $\Sub_G(x)=\langle \bC_G(x),\bN_G(\bS_d)\rangle$ for any $x\in P$.
 \end{enumerate}
\end{prop}

\begin{proof}
Since $\SL_2(q)$ and $\PGL_2(q)$ have non-abelian Sylow 2-subgroups, we must
have $\ell>2$. By \cite[Prop.~2.2]{Ma14}, we have~(b). Hence, any Sylow
$d$-torus has order divisible by the full $\ell$-part of the order of $G$,
giving~(c).  By inspection of the order formulae \cite[Tab.~24.1]{MT11}, (b)
also implies that $\ell$ does not divide the order of the Weyl group of $\bG$
and thus is good for $\bG$ and not a torsion prime (see
\cite[Tab.~14.1]{MT11}), so we have~(a). Finally, as $P$ is characteristic
in $\bS_d^F$, we have $\bN_G(\bS_d)\le\bN_G(\bS_d^F)\le\bN_G(P)$; since
$P$ is abelian, it is contained in a unique Sylow $d$-torus of $\bG$ by
\cite[Prop.~2.2]{CE94} and thus in fact $\bN_G(\bS_d)=\bN_G(P)$. Now~(d) follows
as $\Sub_G(x)=\langle \bC_G(x),\bN_G(P)\rangle$ by \cite[Prop.~2.11]{Ma25}.
\end{proof}

We now dispose of the rank~1 case, that is, when Sylow $\ell$-subgroups of $G$
are cyclic:

\begin{prop}  \label{prop:cyclic}
 Let $d\ge1$ such that $\Phi_d$ divides the generic order of $(\bG,F)$ exactly
 once. Let $\ell{\not|}q$ be a prime with $d=e_\ell(q)$ such that Sylow
 $\ell$-subgroups of $G$ are abelian. Then for all $\ell$-elements
 $1\ne x\in G$, $\Sub_G(x)=\bN_G(\bS_d)=\bN_G(P)$ where $\bS_d$ is a Sylow
 $d$-torus of $\bG$ containing~$x$.
\end{prop}

\begin{proof}
Let $1\ne x\in\bS_d^F$ be an $\ell$-element. Since $\ell$ is good for $\bG$
and not a torsion prime (see
above), $\bC_\bG(x)$ is a $d$-split Levi subgroup of $\bG$ by
\cite[Prop.~2.2]{CE94}. Since $\Phi_d$ divides the generic order of $(\bG,F)$
exactly once, the only $d$-split Levi properly containing $\bC_\bG(\bS_d)$
is~$\bG$ itself. As we assumed $\bG$ to be simple and $\ell$ is not a torsion
prime, we have $x\notin\bZ(G)$, so this forces $\bC_\bG(x)=\bC_\bG(\bS_d)$,
whence we conclude by Proposition~\ref{prop:Burnside}.
\end{proof}

We may and will hence assume in the sequel that $\Phi_d$ divides the order
polynomial of our group at least twice and thus that Sylow $\ell$-subgroups
have rank at least~2 (by \cite[Thm~25.14]{MT11}). In the following two sections
we discuss the exceptional and the classical groups.

\begin{rem}
 Subnormalisers of $\ell$-elements of a given finite permutation group or matrix
 group over a finite field can be determined effectively using the criterion in
 \cite[Cor.~2.10]{Ma25}, which we have implemented in the \GAP\ system
 \cite{GAP}. This will be used throughout to treat small cases.
\end{rem}

\section{Overgroups of Sylow tori normalisers and subnormalisers in groups of exceptional type}   \label{sec:exc}
Throughout this section, $\bG$ is a simple linear algebraic group of simply
connected type and $F:\bG\to\bG$ a Frobenius morphism with respect to an
$\FF_q$-structure. (Subnormalisers in the Suzuki and Ree groups were already
handled in \cite[Thm~5.11]{Ma25}.) We consider primes $\ell$ not dividing $q$.
Since subnormalisers of $\ell$-elements contain the normaliser of a Sylow
$\ell$-subgroup, which in our case, agrees with the normaliser of a Sylow
$d$-torus, we first classify overgroups of the latter type of subgroups.

\begin{thm}   \label{thm:over}
 Assume that $G:=\bG^F$ is of exceptional type, let $d\ge1$ and let
 $\bS_d\le\bG$ be a Sylow $d$-torus. Let $\ell$ be a prime with $d=e_\ell(q)$
 such that Sylow $\ell$-subgroups of $G$ are abelian of rank at least~$2$.
 Then $\bN_G(\bS_d)$ is maximal in the cases listed in Table~\ref{tab:overmax},
 while the proper overgroups of $\bN_G(\bS_d)$ in the other cases are as given
 in Tables~\ref{tab:overgen} and~\ref{tab:overex}.
\end{thm}

The last column in Tables~\ref{tab:overgen} and~\ref{tab:overex}, headed
$\bC_G(x)$, will be explained in Theorem~\ref{thm:subn exc}.

\begin{table}[htb]
\caption{Maximal Sylow $d$-torus normalisers}   \label{tab:overmax}
 $\begin{array}{crc||crc}
   \bG^F& d& \bN_G(\bS_d)& \bG^F& d& \bN_G(\bS_d)\\
   \hline
   \tw3D_4(q)&3,6& \Phi_d^2.G_4& E_8(q)& 1,2& \Phi_d^8.W(E_8)\\
    E_6(q)& 1& \Phi_1^6.W(E_6)& & 3,6& \Phi_d^4.G_{32}\\
       & 3& \Phi_3^3.G_{25}& &   4& \Phi_4^4.G_{31}\\
    \tw2E_6(q)& 2& \Phi_2^6.W(E_6)& & 5,10& \Phi_d^2.G_{16}\\
       & 6& \Phi_6^3.G_{25}& & 12& \Phi_{12}^2.G_{10}\\  
    E_7(q)& 1,2& \Phi_d^7.W(E_7)& \\
 \end{array}$
\end{table}

Here, Sylow tori are indicated by their order polynomial, $W(E_i)$ denotes a
Weyl group of type $E_i$, and $G_i$ denotes a primitive complex reflection
group according to Shephard--Todd.

\begin{table}[htb]
\caption{Generic overgroups of Sylow $d$-torus normalisers}   \label{tab:overgen}
$\begin{array}{ccccc}
 \bG^F& d& \bN_G(\bS_d)& \text{overgroups}& \bC_G(x)\\
\hline\hline
 G_2(q)& 1& \Phi_1^2.W(G_2)& A_2(q).2& \Phi_1.A_1(q)\\
       & 2& \Phi_2^2.W(G_2)& \tw2A_2(q).2& \Phi_2.A_1(q)\\
\hline
 \tw3D_4(q)& 1& \Phi_1^2\Phi_3.W(G_2)& \Phi_3.A_2(q).2& \Phi_1\Phi_3.A_1(q)\\
           & 2& \Phi_2^2\Phi_6.W(G_2)& \Phi_3.\tw2A_2(q).2& \Phi_2\Phi_6.A_1(q)\\
\hline
 F_4(q)& 1,2& \Phi_d^4.W(F_4)& D_4(q).\fS_3& \Phi_d^2.A_2(\pm q)\\
       & 3,6& \Phi_d^2.G_5& \tw3D_4(q).3& \Phi_d.A_2(\pm q)\\
       &   4& \Phi_4^2.G_8& D_4(q).\fS_3& \text{--}\\
\hline
 E_6(q)& 2& \Phi_2^4\Phi_1^2.W(F_4)& \Phi_1^2.D_4(q).\fS_3& \Phi_1^2\Phi_2^2.\tw2A_2(q)\\
       & 4& \Phi_4^2\Phi_1^2.G_8& \Phi_1^2.D_4(q).\fS_3& \text{--}\\
       & 6& \Phi_6^2\Phi_3.G_5& \Phi_3.\tw3D_4(q).3& \Phi_3\Phi_6.\tw2A_2(q)\\
\hline
 \tw2E_6(q)& 1& \Phi_1^4\Phi_2^2.W(F_4)& \Phi_2^2.D_4(q).\fS_3& \Phi_1^2\Phi_2^2.A_2(q)\\
       & 3& \Phi_3^2\Phi_6.G_5& \Phi_6.\tw3D_4(q).3& \Phi_3\Phi_6.A_2(q)\\
       & 4& \Phi_4^2\Phi_2^2.G_8& \Phi_2^2.D_4(q).\fS_3& \text{--}\\
\hline
 E_7(q)& 3,6& \Phi_d^3\Phi_{d/3}.G_{26}& \Phi_{d/3}.E_6(q).2& \Phi_{d/3}\Phi_d.\tw3D_4(q)\\
       &   4& \Phi_4^2.A_1(q)^3.G_8& A_1(q)^3.D_4(q).\fS_3& \text{--}\\
\hline
 E_8(q)&   8& \Phi_8^2.G_9& D_4(q^2).W(G_2)& \Phi_8.A_1(q^4)\\
\end{array}$
\end{table}

\begin{table}[htb]
\caption{Non-generic overgroups of Sylow $d$-torus normalisers}   \label{tab:overex}
$\begin{array}{ccccc}
 \bG^F& d& \bN_G(\bS_d)& \text{overgroups}& \bC_G(x)\\
\hline
 G_2(4)& 2& 5^2.W(G_2)& J_2& 5.A_1(4)\\
 F_4(2)& 4& 5^2.G_8& \tw2F_4(2)& \text{--}\\
 F_4(3)& 6& 7^2.G_5& \tw3D_4(2).3& \text{--}\\
 E_6(2)& 4& 5^2.G_8& \tw2F_4(2),F_4(2)& \text{--}\\
 \tw2E_6(2)& 3& 7^2.3.G_5& F_4(2)\times3& \text{--}\\  
\end{array}$
\end{table}

\begin{proof}
We use knowledge on maximal subgroups of the groups $G$ in question. Note that
$\ell>3$ since Sylow 2- and 3-subgroups of all groups considered are
non-abelian. Thus by our assumptions that $d=e_\ell(q)$, $\ell$ is a Zsigmondy
prime divisor of $\Phi_d(q)$, and $\bN_G(\bS_d)$ contains an abelian Sylow
$\ell$-subgroup of~$G$, by \cite[Thm~25.14]{MT11}, of rank at least~2. The
normalisers of Sylow $d$-tori
are given in \cite[Tab.~1 and~3]{BMM93}. For $G=G_2(q)$ it follows from the
lists in \cite[Tab.~8.30, 8.41 and~8.42]{BHR} (for $p=2$, $p\ge5$ and $p=3$,
respectively) that the only maximal subgroups that can possibly contain
$\bN_G(\bS_d)$ are $\SL_3(\pm q).2$, and in addition $J_2$ for $q=4$ and
$d=2$. By \cite[Tab.~8.3 and 8.5]{BHR} there are no other proper subgroups of
$\SL_3(\pm q).2$ containing $\bN_G(\bS_d)$. Here note that $q\ge8$ for $d=1$
and $q\ne2,3,5,7$ for $d=2$ as otherwise there is no prime $\ell$ as required.
\par
For $G=\tw3D_4(q)$, the result of Kleidman cited in \cite[Tab.~8.51]{BHR} shows
that the only maximal overgroups of $\bN_G(\bS_d)$ for $d=1,2$ are as listed
while for $d=3,6$, $\bN_G(\bS_d)$ is maximal. Again we conclude with
\cite[Tab.~8.3--8.6]{BHR} that the only possible overgroups for $d=1,2$ are the
stated ones.
\par
For $G=F_4(q)$, the results of Craven \cite[Tab.~1, 7 and~8]{Cr23} yield the
given maximal subgroups as possible overgroups. Again by \cite[Tab.~8.50]{BHR},
respectively our previous result for $\tw3D_4(q)$, there are no further
non-maximal overgroups.   \par
For $G=E_6(q)$ the tables \cite[Tab.~2 and~9]{Cr23} allow us to conclude, and
for $G=\tw2E_6(q)$, the tables in \cite[Tab.~3 and~10]{Cr23}. For types $E_7$
and $E_8$ we
use \cite[Thm~0]{Li18} in conjunction with \cite[Thm~1.1 and~1.2]{Cr23a} to
conclude that $\bN_G(\bS_d)$ is maximal in all cases listed in
Table~\ref{tab:overmax}, and to derive the possible overgroups in the other
cases, using also \cite[Thm~29.1]{MT11}. Observe that by assumption these
need to contain an elementary abelian $\ell$-group, with $\ell>5$, of rank at
least~2.
\end{proof}

Let us point out that the generic overgroups in Table~\ref{tab:overgen} are all
obtained as $F$-fixed points of a reductive subgroup $\bH$ of $\bG$ as follows:
the connected component $\bH^\circ$ is generated by the minimal $d$-split Levi
subgroups properly containing $\bN_G(\bS_d)$ corresponding to one conjugacy
class $C$ of reflections in the relative Weyl group $W_d$ of $\bS_d$, and
$(\bH/\bH^\circ)^F\cong W_d/\langle C\rangle$. This is the exact analogue of
the corresponding result for subnormalisers of semisimple elements in simple
algebraic groups in \cite[Thm~6.8]{Ma25}. Nevertheless, the additional
overgroups in Table~\ref{tab:overex} for small values of the parameters
indicate that there may not be any conceptual proof of the preceding
classification, avoiding the precise knowledge of all maximal subgroups.

\begin{thm}   \label{thm:subn exc}
 Let $\bG$ be simple with Frobenius map $F$ with respect to an
 $\FF_q$-structure such that $G=\bG^F$ is of exceptional type. Let $\ell$ be a
 prime dividing $|G|$ such that Sylow $\ell$-subgroups of $G$ are abelian.
 Let $1\ne x\in G$ be an $\ell$-element and $\bS_d\le\bG$ be a Sylow $d$-torus
 containing $x$, where $d=e_\ell(q)$. Then $\Sub_G(x)=G$ unless one of:
 \begin{enumerate}[\rm(1)]
  \item $\bC_G(x)=\bC_G(\bS_d)$, $x$ is picky and $\Sub_G(x)=\bN_G(\bS_d)$; or
  \item $\bC_G(x)>\bC_G(\bS_d)$ is the $d$-split Levi subgroup in the last
   column of Tables~\ref{tab:overgen} or~\ref{tab:overex}, and $\Sub_G(x)$ is
   as in the next to last column.
 \end{enumerate}
\end{thm}

\begin{proof}
Since Sylow $\ell$-subgroups of $G$ are supposed to be abelian, we have
$\Sub_G(x)=\langle\bC_G(x),\bN_G(P)\rangle$ by Proposition~\ref{prop:Burnside}.
By \cite[Thm~5.3]{Ma25} the element $x$ is picky if and only if we are in
Case~(1), and then $\Sub_G(x)$ is as stated. Now assume $\bC_G(x)$ properly
contains $\bC_G(\bS_d)$. Since $\ell$ is a good prime for $\bG$, then
$\bC_\bG(x)$ is a $d$-split Levi subgroup of $\bG$, see \cite[Prop.~2.2]{CE94},
properly containing $\bC_\bG(\bS_d)$. The possibilities can be obtained in
\Chevie\ \cite{Mi15} and are also listed in \cite[Tab.~3.3]{GM20} (unfortunately
with some omissions for $E_8$ when $d=8,12$). By \cite[Prop.~5.2]{Ma25},
$\Sub_G(x)$ is generated by the $F$-fixed points of the normalisers of the
Sylow $d$-tori of $\bG$ containing $x$. Thus, either $\Sub_G(x)=G$, or it is
one of the groups in Theorem~\ref{thm:over}. In the four cases with an entry
``--" in the last column in Table~\ref{tab:overgen}, the stated overgroup does
not contain any $d$-split Levi subgroup of $\bG$ properly containing
$\bN_G(\bS_d)$, so here $\Sub_G(x)=G$.
In all other cases, the listed $d$-split Levi is the unique one embedding in
the listed overgroup. In Table~\ref{tab:overex} it is easy to see from the
known character tables and information on maximal subgroups that only in the
stated case the listed overgroup of $\bN_G(\bS_d)$ does contain the
centraliser of an $\ell$-element $x$ with $\bC_G(x)>\bC_G(\bS_d)$. Note that
the prime $\ell$ is uniquely determined in each of the cases in
Table~\ref{tab:overex}, namely $\ell=5,5,7,5,7$ in the respective cases.
In $G_2(4)$, a \GAP-computation shows that both the generic case $\SU_3(4).2$
and the exotic case $J_2$ do occur as subnormalisers for suitable 5-elements.
\end{proof}

We can check Conjecture~\ref{conj:AN} for the exceptional subnormaliser:

\begin{prop}
 Conjecture~\ref{conj:AN} holds for $G_2(4)$ at $\ell=5$ with $\Sub_G(x)=J_2$.
\end{prop}

\begin{proof}
From the known character tables one sees
that both $G$ and $\Sub_G(x)=J_2$ possess 14 irreducible character of degree
prime to~5 not vanishing on $x$, and four further ones of degree divisible
by~5 exactly once, and there is a bijection such that the values at~$x$ of
corresponding characters agree up to sign, so in particular generate the same
field extension.
\end{proof}

\section{Overgroups of Sylow tori normalisers and subnormalisers in groups of classical type}   \label{sec:class}
Here, throughout we let $\bG$ be simple of simply connected classical type and
$F:\bG\to\bG$ a Frobenius map with respect to an $\FF_q$-structure, not
inducing triality. Let $G:=\bG^F$. Let $\ell>2$ be a prime not dividing $q$
and $H$ an elementary abelian $\ell$-subgroup of $G$ of maximal
possible rank. As $\ell$ is not a torsion prime for $\bG$, $H$ embeds into
a maximal torus $\bT$ of~$\bG$ by \cite[Cor.~14.17]{MT11}. Thus, $\bT$ lies
in the $F$-stable subgroup $\bC_\bG(H)$, which possesses $F$-stable maximal
tori, all of which contain $H$. Hence, without loss we may assume $\bT$ to be
$F$-stable.
Furthermore, $\bT$ then contains a Sylow $d$-torus $\bS_d$ of $(\bG,F)$ for
$d=e_\ell(q)$, by \cite[Prop.~5.7]{Ma07}, and $H\le\bS_d$. Thus $H$ is
normalised by $N:=\bN_G(\bS_d)$. We are thus led to studying overgroups of
$\bN_G(\bS_d)$. This we will do in terms of Aschbacher's classification, by
which any maximal subgroup either lies in a collection $\cC_i$ of natural,
geometric subgroups, or else is almost simple modulo its centre and lies in a
class denoted~$\cS$ (see \cite{KL90}, \cite{BHR}, or \cite[\S27]{MT11} for an
introduction and further references).

In all types we denote by $P_r$ a maximal parabolic subgroup of $G$ stabilising
an $r$-dimensional totally singular subspace of the natural $\FF_q$-
respectively $\FF_{q^2}$-module.

\subsection{The special linear groups}   \label{subsec:SL}

We start our investigation with the groups $\SL_n(q)$. We fix the following
setup and notation throughout this subsection. Let $G=\SL_n(q)$ with $n\ge2$,
$(n,q)\ne(2,2),(2,3)$. Let $d\ge1$ be an integer and write $n=ad+r$ with
$0\le r<d$. Let $\bS_d<\bG=\SL_n$ be a Sylow $d$-torus. To avoid certain
degenerate situations and since this will be satisfied for our intended
application, we also assume that there is some Zsigmondy primitive prime
divisor $\ell>2$ of $q^d-1$. So, in particular, $\ell$ divides $|G|$ if
$d\le n$. We first consider maximal subgroups of $G$ of geometric type (in
Aschbacher's approach).

\begin{prop}   \label{prop:over SL}
 Let $G=\SL_n(q)$, and $d,\ell$ as above. If $M<G$ is a maximal subgroup
 containing the normaliser $N:=\bN_G(\bS_d)$ of a Sylow $d$-torus
 $\bS_d\le\bG$, then one of:
 \begin{enumerate}[\rm(1)]
  \item $M=P_{ad}$ or $P_r$ is maximal parabolic, if $r>0$; 
  \item $M=(\GL_d(q)\wr\fS_a)\cap G$ if $n=ad$ and $a>1$;
  \item $M=\GL_{n/t}(q^t).t\cap G$ if $n=d$ and $t|n$ is a prime;
  \item $M=\Sp_n(2)$ if $n>2$ is even and $q=d=2$;
  \item $M=\SU_n(2)$ if $n\ge3$, $q=4$ and $d=1$; or
  \item $M$ is in class $\cS$.
 \end{enumerate}
\end{prop}

Note that in cases (4) and (5) we necessarily have $\ell=3$.

\begin{proof}
By \cite[Exmp.~3.5.14]{GM20} the centraliser and normaliser of $\bS_d$ have
the structure
$$\bC_G(\bS_d)=\big(\GL_1(q^d)^a\times\GL_r(q)\big)\cap G\quad\text{and}\quad
  N=\big(\GL_1(q^d)\wr G(d,1,a)\times\GL_r(q)\big)\cap G$$
(all viewed as subgroups of $\GL_n(q)$).
Observe that for $d=1$, $N$ contains a subgroup $N_0:=C_\ell^{n-1}\rtimes\fS_n$,
where the base group is the deleted permutation module for the symmetric group
$\fS_n$, while for $d>1$ it contains a subgroup $N_0:=C_\ell^a\rtimes G(d,1,a)$
where here the complement acts irreducibly on the base group. In either case,
the $N$-composition factors of the natural module have dimensions
$ad$ and $r$, and $N$ acts primitively on the factor of dimension~$r$.
For the proof we now go through the various Aschbacher classes of maximal
subgroups of $G$ as described, for example, in \cite[Prop.~28.1]{MT11}.

Assume first $M=P_m$, with $1\le m<n$, is a maximal parabolic subgroup with
$N\le M$, so $M$ acts maximally reducibly. Then by what we just observed,
$m\in\{ad,r\}$, and $P_m$ is proper if $r>0$, so we arrive at case~(1).

Next assume $M=(\GL_m(q)\wr\fS_t)\cap G$, with $n=mt$, $t\ge2$, is imprimitive.
If $d=1$ then the automiser in $M$ of an elementary abelian $\ell$-subgroup $E$
of $\GL_m(q)^t\cap G$ of rank~$n-1$ is $\fS_m\wr\fS_t$, a proper subgroup of
the automiser $\fS_n$ in $G$, unless $m=1$ as in case~(2). If $d>1$ write
$m=a_1d+r_1$ with $0\le r_1<d$. Then $M$ has $\ell$-rank~$ta_1$, so we need
$ta_1\ge a$ whence $ta_1=a$. Then the automiser in $M$ of $E$ is
$G(d,1,a_1)\wr\fS_t$, again a proper subgroup of $G(d,1,a)$ unless $a_1=1$, so
$t=a$ and hence $r=ar_1$. For its centraliser in $M$ to contain a subgroup
$\GL_r(q)$ we need $r=0$ and so again obtain~(2).

If $M=\GL_m(q^t).t\cap G$, with $n=mt$, $t$ prime, is an extension field
subgroup, then $M$ contains an elementary abelian $\ell$-group of rank~$a$,
respectively of rank~$n-1$ if $d=1$, only if $t|\gcd(r,d)$. Its automiser in
$M$ is then $G(d/t,1,a).t$, which equals $G(d,1,a)$ only if $a=1$. Comparing
the centralisers of a Sylow $d$-torus in $M$ and $G$ we then see that $r=0$.
So we must have $n=d$, giving~(3).

Next assume $M=\GL_{n_1}(q)\otimes\GL_{n_2}(q)\cap G$ with $n_1n_2=n$,
$1<n_1< n_2$, preserves a tensor decomposition. Writing $n_i=a_id+r_i$ with
$0\le r_i<d$, the $\ell$-rank of $M$ is $a_1+a_2$, while $G$ has $\ell$-rank
at least $a_1a_2d+a_1r_2+a_2r_1-1$, and this is bigger unless $r_1=r_2=0$ and
either $a_1=a_2=1$, $d=2$, so $n_1=n_2=2$ which is excluded, or $a_1=2,a_2=3$,
$d=1$. In the latter case the automiser in $G$ is $\fS_6$, while it is only
$\fS_2\times\fS_3$ in $M$. The same reasoning applies to $M=\GL_m(q)\wr\fS_t$
with $n=m^t$, $t\ge2$, $m\ge3$.

Now let $M=t^{2m+1}.\Sp_{2m}(t)\bZ(G)$ be the normaliser of an extra-special
$t$-subgroup with $n=t^m$, where furthermore $q$ is minimal among powers of $p$
with $q\equiv1\pmod {t(2,t)}$. The maximal order of an abelian subgroup in
$M$ is $t^{m(m+1)/2+m+1}\gcd(n,q-1)$ by \cite[Thm~3.1]{Vd99}, while $N$
contains a maximal torus of order at least $(q-1)^{n-1}$. The ensuing
inequality does not hold for any $n\ge5$.

For $M=\Sp_n(q)$ with $n$ even to have large enough $\ell$-rank, the order
formulas show that we need $d$ to be even. In this case, $M$ contains a
subgroup isomorphic to $N_0$, in the normaliser of a Sylow $d$-torus of $M$.
But the semisimple part of its centraliser in $M$ is a symplectic group
$\Sp_r(q)$, properly smaller than $\GL_r(q)$ unless $r=0$ (note that here $r$
is necessarily even as both $n,d$ are). In the latter case,
$\bC_G(\bS_d)$ contains a subgroup of index $q-1$ of a homocyclic group
$(q^d-1)^a$, with $n=ad$, and such an abelian subgroup can only be contained
in~$M$ if $q^{d/2}-1=1$, so $q=d=2$ (and thus $\ell=3$) as in~(4).

If $M=\SO_n^{(\pm)}(q)$ with $q$ odd, arguing as in the previous case we reach
the contradiction $q=2$. For $M=\SU_n(q_0)$ with $q=q_0^2$ the same type of
consideration leads to possibility~(5).

Finally, assume $M=\GL_n(q_0)\cap G$, with $q=q_0^f$, $f\ge2$, is a subfield
group. Note that $N$ contains a maximal torus of $G$, of order at least
$(q^d-1)^a(q-1)^{n-da-1}$, while the size of maximal tori in $M$ is bounded
above by $(q_0+1)^{n-1}$, with $q_0^2\le q$. The resulting inequality forces
$q=4$, $q_0=2$ and thus $\ell=3$, $d=1$, but in this case the $3$-rank of $M$
is at most $\lfloor (n+1)/2\rfloor$ while $N_0$ has $3$-rank~$n-1$, forcing
$n\le3$. But for $n=2$, $M$ acts imprimitively and already occurs under~(2),
while for $n=3$ it does not contain a subgroup $N_0$.
\end{proof}

In order to deal with the maximal subgroups in class $\cS$, we first derive
some easy estimates.

\begin{lem}   \label{lem:|SL_n|}
 Let $n\ge2$ and $q$ be a prime power. Then
 $|\SL_n(q)|>(q-1) q^{n^2-2}$.
\end{lem}

\begin{proof}
By Euler's pentagonal number theorem, for all real $x$ with $|x|<1$ we have
$$\prod_{i=1}^\infty (1-x^i)
  =1+\sum_{k=1}^\infty (-1)^k\big(x^{k(3k+1)/2}+x^{k(3k-1)/2}\big).$$
In particular, if $0\le x\le 1/2$ then
$$\prod_{i=1}^\infty (1-x^i)
  =1-x-x^2+x^5+x^7\ldots\ge 1-x-x^2+x^5=(1-x)(1-x^2-x^3-x^4).$$
Thus, for $x:=q^{-1}$ and $n\ge2$ we obtain
$$\prod_{i=2}^n\frac{q^i-1}{q^i}=\prod_{i=2}^n(1-q^{-i})\ge
  1-q^{-2}-q^{-3}-q^{-4} > \frac{q-1}{q},$$
showing that
$$\prod_{i=2}^n(q^i-1) > (q-1)q^{\binom{n+1}{2}-2}.$$
Our claim is immediate from this with the order formula for $\SL_n(q)$.
\end{proof}

\begin{lem}   \label{lem:N in SL}
 The subgroup $N$ of $G=\SL_n(q)$ has order $|N|\ge q^{n-1}$.
\end{lem}

\begin{proof}
By our assumptions we have $d>1$ if $q\le3$, and thus $q^d-1\ge 3q^d/4$.
Then since $|G(d,1,a)|=d^a|\fS_a|=d^a a!\ge (4/3)^a$ we conclude that
$(q^d-1)^a|G(d,1,a)|\ge q^{ad}$ unless $ad=1$. But $ad>1$ in our
situation. Since $|N|=(q^d-1)^a|G(d,1,a)|\cdot|\SL_r(q)|$, our claim follows
from Lemma~\ref{lem:|SL_n|} when $r\ge2$. Direct computation shows that it
also holds when $r\in\{0,1\}$.
\end{proof}

\begin{prop}   \label{prop:S in SL}
 In the situation of Proposition~\ref{prop:over SL}, assume $M$ is a maximal
 subgroup of $G$ in class $\cS$ containing $N$. Then $M<G$ is one of:
 \begin{enumerate}[\rm(1)]
  \item $2.\fA_5<\SL_2(9)$, with $d=2$, $\ell=5$;
  \item $2.\fA_5<\SL_2(11)$, with $d=1$, $\ell=5$;
  \item $3.\fA_6<\SL_3(4)$, with $d=1$, $\ell=3$ or $d=2$, $\ell=5$; or
  \item $\fA_7<\SL_4(2)$, with $d=3$, $\ell=7$.
 \end{enumerate}
\end{prop}

\begin{proof}
We consider the various possibilities for the non-abelian simple composition
factor $S:=F^*(M)/\bZ(F^*(M))$ of $M$ according to the classification. First
assume $S$ is not of Lie type in the same characteristic as $G$.

(1) For $n\le12$ the possible $M\in\cS$ are given in the tables of
\cite[\S8]{BHR}. By Lemma~\ref{lem:N in SL}, if $N\le M$ we must have
$$|\Aut(S)|\ge |N|/|\bZ(G)|\ge q^{n-1}/|\bZ(G)|\ge q^{n-1}/n,$$
which for any $S$ gives a small
upper bound on~$q$. Also, $\bar N:=N\bZ(G)/\bZ(G)$ contains elements of order
at least $q-1$, so $\Aut(S)$ must contain elements of at least that order,
yielding further restrictions on~$q$. The occurring subgroups $M\in\cS$ can
now be investigated using the Atlas \cite{Atl}, leading to the four items
listed in the conclusion.
Hence from now on we can assume $n\ge13$. Then $N$ contains a maximal torus of
$G$ of size at least $(q^2-1)^6\ge3^6$ when $q\le3$, respectively
$(q-1)^{11}\ge 3^{11}$ when $q\ge4$,

(2) Now assume $S=\fA_m$ is an alternating group. If $S$ has a faithful
projective representation of degree less than $m-2$ then $m\le 8$ by
\cite[Prop.~5.3.7]{KL90} so since we have $n\ge13$, we may assume $m\le n+2$.
As argued in the proof of Lemma~\ref{lem:N in SL} then $q^d-1\ge 3q^d/4$.
If $r\ge2$ then $\bar N$ contains an abelian subgroup of size~$(q^d-1)^a$,
and $(q^d-1)^a\ge(4^d-1)^a\ge (4-1)^{ad}=3^{ad}$ if $q\ge4$.
But the maximal size of an abelian subgroup of $\fS_m$ is bounded above by
$3^{m/3}$ by \cite[Thm~1]{BG89}, and since $ad>n/2$ we conclude we must have
$$3^{n/2}< 3^{ad} < 3^{(n+2)/3}$$
whence $n<4$, a contradiction. If $q=3$ the existence of $\ell$ forces
$d\ge3$ and then $3^d-1\ge 3^{5d/6}$. So in this case we deduce
$(3^d-1)^a\ge 3^{5ad/6}$ must be less than $3^{5n/12}\le 3^{(n+2)/3}$, whence
$n<8$. If $r\le1$ and still $q\ge3$ then $\bar N$ will still contain an abelian
subgroup of size at least $(q^d-1)^a/n(q-1)$, but we also have $ad\ge n-1$ and
again the required inequality is not satisfied for $n>11$.

It remains to consider the possibility that $q=2$. First assume $r\le2$.
As $d\ge2$ we have $2^d-1\ge 2^{4d/5}$ and so conclude we must have
$2^{4(n-2)/5}\le 3^{(n+2)/3}$, so $4(n-2)/5<8(n+2)/15$ which implies $n<10$,
a case already considered. Finally, if $r>2$ then we use that
$N$ contains an elementary abelian $\ell$-group of rank~$a$ centralised by
a subgroup $\SL_r(2)$. Now the centraliser in $\fS_m$ of such an elementary
abelian $\ell$-subgroup is $\fS_{m-a\ell}$ times an $\ell$-group, so
$\SL_r(2)$ must embed into $\fS_{m-a\ell}$ and hence possess a faithful
representation in characteristic~0 of degree at most $m-a\ell-1$. This implies 
$2^{r-1}\le m-a\ell-1$ by \cite{TZ96}. On the other hand, as $\ell$ is a
Zsigmondy prime divisor of $q^d-1$ we have $\ell\ge d+1$. So
$$2^{r-1}\le m-a\ell-1\le n+2-a(d+1)-1=1-a+r\le r$$
which is never satisfied for $r>2$.

(3) Next assume $S$ is of Lie type in characteristic not dividing~$q$. For
$S=\PSL_m(y)$ with $m\ge3$ the
smallest degree of a faithful projective representation in characteristic
not dividing $y$ is $(y^m-1)/(y-1)-m\ge y^{m-1}$ by \cite[Tab.~1]{TZ96},
unless $(m,y)\in \cE:=\{(3,2),(3,4),(4,2),(4,3)\}$. Furthermore,
$$|M|\le|\Aut(S)|\le 2|\PGL_m(y)|\log_2 y\le 2\ y^{m^2-1}\log_2 y.$$
Thus, for $(m,y)\notin \cE$, if $N\le M$ then
$$q^{y^{m-1}-1}\le q^{n-1}\le |N|\le |M|\le 2 y^{m^2-1}\log_2 y$$
by Lemma~\ref{lem:N in SL}. This is only satisfied for $(m,y)=(3,3)$. It thus
remains to discuss $(m,y)\in \cE\cup\{(3,3)\}$. All projective irreducible
representations of $\PSL_3(2)$ have degree at most~8, the group $\PSL_3(3)$
does not have an irreducible representation of degree $n>12$ satisfying the
inequality, for $\PSL_3(4)$ and $\PSL_4(2)\cong\fA_8$ the inequality only
holds when $n\le12$, and for $\PSL_3(4)$ it is never satisfied. Thus no
further cases arise.

For $S=\PSL_2(y)$, $y\ge7$ and $y\ne9$ (as we already considered alternating
groups), the normaliser of a (cyclic) Sylow $\ell$-subgroup has order at
most $y(y-1)$, while the minimal projective degree is $(y-1)/(2,y-1)$, so we
need
$$y(y-1)\ge q^{n-1}\ge q^{(y-1)/2}-1.$$
This is not satisfied for any $n\ge13$ and prime powers $y,q$, so no further
examples arise.

For $S=\PSU_m(y)$ with $m\ge3$ the smallest projective faithful degree is at
least $(y-1)y^{m-2}$, unless $(m,y)\in\{(4,2),(4,3)\}$, by \cite[Tab.~I]{TZ96},
and $|\Aut(S)|\le 2y^{m^2-1}\log y$. Arguing as above, the only cases that
satisfy the relevant inequality are
$$(m,y)\in\{(3,3),(3,4),(3,5),(5,2),(6,2)\}.$$
A maximal abelian subgroup of $\PSU_3(y)$ has size $(y+1)^2$ by
\cite[Thm~3.1]{Vd99}. But as pointed out above, for $n\ge13$ there is a torus
in $N$ of size at
least $(q^2-1)^6$ when $q\le3$, respectively $3^{11}$ when $q\ge4$, hence the
case $m=3$ does not lead to new examples. For $S=\PSU_4(2)$ or $\PSU_4(3)$,
$\Aut(S)$ contains no abelian subgroups of the required order.
For $S=\PSU_5(2)$ the relevant inequality only holds when $q=3$ and $n\le16$,
but by \cite{HM} the smallest degree $n>12$ of a 3-modular projective
irreducible representation of~$S$ is~44. For $S=\PSU_6(2)$ the inequality only
holds when $q=3$ and $n\le23$, while by \cite{HM} the only 3-modular projective
irreducible representation of $S$ in degree $12<n\le23$
has degree~21, in which case $|\bS_d^F|$ is too large. 

For $S=\PSp_{2m}(y)$ with $m\ge2$ the smallest faithful projective degree is
$(y^m-1)/2$ if $y$ is odd, and $(y^m-1)(y^m-y)/(2(y+1))$ if $y$ is even.
(Note that we may assume $(m,y)\ne(2,2)$ as $\PSp_4(2)\cong\fS_6$.)
Our inequality (with $n>12$) is satisfied only for
$(2m,y)\in\{(4,3),(4,5),(4,7),(6,2),(6,3),(8,3)\}$. The condition that
$\Aut(S)$ should contain an abelian subgroup of size $|\bS_d^F|/\bZ(G)$ rules
out all but $S=\PSp_6(3)$ with $q=2$. The smallest faithful irreducible
2-modular representations of $S$ have dimension~13 and~78. But those of
degree~13 are only defined over $\FF_4$, while the order of $N$ in
$\SL_{78}(2)$ is much too large.

For $S=\PSO_{2m+1}(y)$ with $m\ge3$ and $y$ odd by \cite[Tab.~1]{TZ96} the
smallest projective degree is at least $(y^m-1)(y^m-y)/(y^2-1)$, unless
$(m,y)=(3,3)$, and the necessary inequality is never satisfied. For
$S=\PSO_7(3)$ the smallest faithful projective degree is~27 by \cite{HM},
but $\Aut(S)$ has no abelian subgroup of size at least $3^{12}$.

For $S=\PSO_{2m}^+(y)$ with $m\ge4$, again by \cite{TZ96} the minimal
projective degree is at least $(y^m-1)(y^{m-1}-1)/(y^2-1)-7$, respectively~8
for $\PSO_8^+(2)$. Our inequality is never satisfied in the former case; for
$S=\PSO_8^+(2)$ we obtain $q=3$, but by \cite{HM} there is no projective
irreducible 3-modular representation of $S$ of degree $n\ge13$ for which the
inequality would be satisfied.

For $\PSO_{2m}^-(q)$ with $m\ge4$ the minimal degree is at least
$(y^m+1)(y^{m-1}-y)/(y^2-1)-m+2$ and our inequality is never satisfied.
For $S$ of exceptional Lie type, the lower bounds on faithful projective
representations in \cite{TZ96} are always large enough to exclude these
possibilities.

(4) Now assume $S$ is sporadic. Since $N$ and thus $M$ contain elements of
order $(q^d-1)/(q-1)$ we see that $d\le6$, and $d\le5$ if $q\ge3$. For $q=2$,
$d=6$, there is no Zsigmondy prime, so we have $d\le5$ and hence $a\ge2$ as
$n\ge13$. Thus $\ell^2$ divides $|S|$, forcing $\ell\le13$. By inspection, if
$\ell^a$ divides $|\Aut(S)|$ then $n<(a+1)d\le(a+1)(\ell-1)\le50$. Comparing
to the list of degrees of faithful irreducible projective representations of
$S$ below $50$ in \cite{HM} shows that no examples arise.

(5) Finally assume $S$ is of Lie type in the same characteristic as~$G$. Again,
when $n\le12$ the tables in \cite[\S8]{BHR} show that no example exists,
so assume $n\ge13$. If $S=\PSL_m(p^f)$ then $n\ge m(m-1)/2$ and $p^f|q$,
or $n\ge m^k$ and $p^f|q^k$ for some $k\ge2$ by \cite[Prop.~5.4.6,
5.4.11]{KL90}. Now write $m=bd+s$ with $0\le s<d$, so that $S$ has $\ell$-rank
at most $b$. If $d=1$ then $G$ has $\ell$-rank at least $n-1\ge m(m-1)/2-1>m=b$,
a contradiction, so $d\ge2$. Then
$n\ge m(m-1)/2=(bd+s)(bd+s-1)/2\ge d(b^2d/2+bs-b/2)$, so $G$ has $\ell$-rank at
least $b^2d/2+bs-b/2>b$, again a contradiction unless $d=2,b=1,s=0$, but then
$n\le12$.

If $S=\PSp_{2m}(p^f)$ then $n\ge 2m(2m-1)-2$ and $p^f|q$, or $n\ge(2m)^k$ and
$p^f|q^k$ for some $k\ge2$, or $m\le6$ and $n=2^m$, by \cite[Tab.~5.4.A]{KL90}.
The first two cases are treated as before. To exclude $n=2^m$, with
$4\le m\le 6$, observe that the $\ell$-rank of $S$ is at most $m/(d/2)$, and
that of $G$ is at least $2^m/d-1>2m/d$.
The same line of argument now applies to all types of groups $S$,
using the bounds in \cite[Tab.~5.4.A, 5.4.B]{KL90}. For $S$ a triality group
or of (possibly twisted) type $B_2,G_2$ or $F_4$ we also refer to the
description in \cite[Rem.~5.4.7]{KL90} for fields of definition. No further
examples arise.
\end{proof}

\begin{rem}   \label{rem:uniq SL}
 The previous results show, in particular, that in the situation of
 Proposition~\ref{prop:over SL} the normaliser of a Sylow $d$-torus lies in a
 \emph{unique} maximal subgroup of $\SL_n(q)$ whenever Sylow $\ell$-subgroups
 are non-cyclic and $r=0$, and $\ell\ne3$ when $q\in\{2,4\}$. Note that
 in the latter exceptions, Sylow 3-subgroups of $\SL_n(q)$ are non-abelian.
 (Compare to the cases in \cite[Tab.~B]{BBGT}.)
\end{rem}

We can now determine the subnormalisers:

\begin{thm}   \label{thm:subn SL}
 Let $G=\SL_n(q)$ with $n\ge2$, let $\ell\nmid q$ be a prime such that Sylow
 $\ell$-subgroups of $G$ are abelian and $\bS_d\le\bG$ a Sylow $d$-torus,
 where $d=e_\ell(q)$. Write $n=ad+r$ with $0\le r<d$. Then for $x\in\bS_d^F$
 an $\ell$-element we have
 \begin{enumerate}[\rm(1)]
  \item $\Sub_G(x)=\bN_G(\bS_d)$ if $\bC_G(x)=\bC_G(\bS_d)$; or
  \item $\Sub_G(x)=\big(\GL_{ad}(q)\times\GL_r(q)\big)\cap G$ if $r>0$ and
   $\bC_{\GL_n(q)}(x)=\prod_i \GL_{n_i}(q^d)\times\GL_r(q)$ with $a=\sum n_i$
   and at least one $n_i>1$; or
  \item $\Sub_G(x)=\big(\GL_d(q)\wr\fS_a\big)\cap G$ if $r=0$, $d>1$, $a>1$ and
   $\bC_{\GL_n(q)}(x)= \GL_1(q^d)^{a-1}\times\GL_d(q)$; or
  \item $\Sub_G(x)=G$ otherwise.
 \end{enumerate}
\end{thm}

\begin{proof}
By Proposition~\ref{prop:Burnside} the subnormaliser of~$x$ is generated by
$N:=\bN_G(\bS_d)$ and $\bC_G(x)$. Since $\bG=\SL_n$ is simply connected and
Sylow $\ell$-subgroups of $G$ are abelian, the centraliser $\bC_\bG(x)$ is a
$d$-split Levi subgroup of $(\bG,F)$ by \cite[Prop.~2.2]{CE94}. First note
that if $d>n/2$ then $\Phi_d$ divides the order polynomial of $\SL_n$ just
once, and so $\Sub_G(x)=\bN_G(\bS_d)$ and $\bC_G(x)=\bC_G(\bS_d)$ by
Proposition~\ref{prop:cyclic} unless $x=1$, so we reach~(1) or~(4) of the
conclusion. The same holds if $\bC_G(x)=\bC_G(\bS_d)\le\bN_G(\bS_d)$ is
a minimal $d$-split Levi.
\par
If $d\le n/2$ we go through the possible maximal overgroups of $\bN_G(\bS_d)$
classified in Propositions~\ref{prop:over SL} and~\ref{prop:S in SL} to see
which ones can possibly contain a non-minimal $d$-split Levi subgroup $L$
of~$G$. By \cite[Exmp.~3.5.14]{GM20} the latter have the form
$$\Big(\GL_{n_1}(q^d)\times\cdots\times\GL_{n_t}(q^d)\times\GL_s(q)\Big)
  \cap G\qquad\text{with}\quad d\sum_i n_i+s=n$$
(which of course implies $s\ge r$). As $L$ is non-minimal, we may assume that
$s>r$ or $n_1>1$, say. First consider the overgroups in class $\cS$ in
Proposition~\ref{prop:S in SL}. Only $M=3.\fA_6<\SL_3(4)$ with $d=1$, $\ell=3$
has non-cyclic Sylow $\ell$-subgroups. By explicit computation in \GAP, this
does not occur as a subnormaliser of a 3-element in $\SL_3(4)$.

Now consider the geometric subgroups of $G$ classified in
Proposition~\ref{prop:over SL}, and first assume $r=0$. Then Cases~(2)--(5)
from Proposition~\ref{prop:over SL} are relevant. In Case~(3) we have $d=n$,
so Sylow $\ell$-subgroups of $G$ are cyclic, a situation already discussed
before. In Cases~(4) and~(5) we have $\ell=3$ and Sylow $\ell$-subgroups of~$M$
are abelian and non-cyclic only for $M=\Sp_4(2)<\SL_4(2)\cong\fA_8$, but direct
computation shows that $M$ is not the subnormaliser of any 3-element of
$\SL_4(2)$. Finally, in Case~(2) the action on the natural module shows that
for $\bC_G(x)$ to be contained in $M$ we need the structure given in (3) of the
conclusion. Here, since $\Sub_G(x)$ contains the $\GL_d(q)$-factor from
$\bC_G(x)$ as well as the symmetric group $\fS_a$ from the normaliser of a
Sylow $d$-torus, we see that $\Sub_G(x)=M$, as claimed. Note that if $d=1$ we
are in Case~(1) and if $a=1$ we are in Case~(4).

So, finally assume $r>0$. Then only the parabolic subgroups
$M\in\{P_{ad},P_r\}$ occur in Proposition~\ref{prop:over SL}. Comparing
dimensions of composition factors on the
natural module we see that either can contain $L=\bC_G(x)$ as above only if
$r=s$. Now note that the $\GL_r(q)$-factor times its centraliser in $G$ is
conjugate to a Levi factor of $M$, and it contains the normaliser of a Sylow
$d$-torus. Thus, if $\Sub_G(x)\le M$ then it already lies in a Levi factor.
Thus $\Sub_G(x)$ has the form $(H\times\GL_r(q))\cap G$ for some subgroup $H$
of $\GL_{ad}(q)$ containing a Sylow $d$-torus normaliser. So we can again
appeal to Proposition~\ref{prop:over SL} to see that if $H$ is proper, it must
lie in a subgroup of type~(2)--(6). Here note that the subgroups of type~(2)
act imprimitively on a sum of $a$ subspaces of dimension $d$, while $L$ has a
primitive summand of dimension at least $2d$ (as $n_1>1$), so this case is out.
In Case~(3) the Sylow $\ell$-subgroups are cyclic, contrary to assumption, and
Cases~(4)--(6) do not occur by the same arguments as given in the previous
paragraph.
So here $\Sub_G(x)$ is a Levi factor of $M$, as in~(1) of the conclusion.
\end{proof}

\subsection{The special unitary groups}

We next consider the special unitary groups. So throughout this subsection let
$G=\SU_n(q)$ with $n\ge3$, $(n,q)\ne(3,2)$. Let $d\ge1$ be an integer and
$$e:=\begin{cases} 2d& \text{if $d$ is odd},\\
              d/2& \text{if $d\equiv2\pmod4$},\\
              d& \text{if $d\equiv0\pmod4$.}\end{cases}$$
Write $n=ae+r$ with $0\le r<e$. We again assume there is a Zsigmondy primitive
prime divisor $\ell>2$ of $q^d-1$; so we have $d=e_\ell(q)$ and $e=e_\ell(-q)$.
Let $\bS_d\le\bG$ be a Sylow $d$-torus with normaliser $N:=\bN_G(\bS_d)$.
By \cite[Exmp.~3.5.14]{GM20} we have
$$\bC_G(\bS_d)=\big(\GL_1((-q)^e)^a\times\GU_r(q)\big)\cap G\quad\text{and}\quad
  N=\big(\GL_1((-q)^e)\wr G(e,1,a)\times\GU_r(q)\big)\cap G.$$

\begin{prop}   \label{prop:over SU}
 Let $G=\SU_n(q)$ with $n\ge3$ and $d,e,\ell$ as above. If $M<G$ is a maximal
 subgroup containing the normaliser $N:=\bN_G(\bS_d)$ of a Sylow $d$-torus
 $\bS_d\le\bG$, then one of:
 \begin{enumerate}[\rm(1)]
  \item $M=(\GU_{ae}(q)\times\GU_r(q))\cap G$ if $r>0$;
  \item $M=(\GU_e(q)\wr\fS_a)\cap G$ if $n=ae$ and $a>1$;
  \item $M=\GL_{n/2}(q^2).2\cap G$ if $n=e$ is even;
  \item $M=\GU_{n/t}(q^t).t\cap G$ if $n=e$ and $2<t|n$ is a prime;
  \item $M=\Sp_4(q).(q-1,2)$ if $n=d=4$ and $q\le3$; or
  \item $M$ is in class $\cS$.
 \end{enumerate}
\end{prop}

\begin{proof}
The situation for $\SU_n(q)$ is Ennola dual to the one for $\SL_n(q)$.
The arguments are now similar to the case of $\SL_n(q)$ in
Proposition~\ref{prop:over SL}, where as far as divisibility questions are
concerned, we need to replace $q$ by $-q$ and $d$ by~$e$, and we now appeal to
\cite[Tab.~3.5B]{KL90} for the description of the Aschbacher classes.

Assume first $M$ is reducible. Since the composition factors of $N$ on the
natural module have dimensions $ae$ and $r$, we either have $M$ is as in~(1),
or $M$ is a parabolic subgroup $P_r$ with $r>0$ (note that $ae>n/2$ is larger
than the dimension of a totally isotropic subspace). But the latter has three
composition factors on the natural module for $G$.

The argument for the imprimitive groups $M=(\GU_m(q)\wr\fS_t)\cap G$ with
$n=mt$ is identical to the one for $\SL_n(q)$ and we reach conclusion~(2).
The subgroups $\GL_{n/2}(q^2).2\cap G$ are normalisers of Levi subgroups
of~$G$ and it can be seen from the description of their Sylow $d$-normaliser
in the proof of Proposition~\ref{prop:over SL} that the automiser of a maximal
rank $\ell$-subgroup is strictly smaller than $G(e,1,a)$ unless $a=1$, so $n=e$
as in~(3).

By arguments as in the case of $\SL_n(q)$, the only further extension field
subgroups that can occur are as in~(4), while again there are no examples for
stabilisers of tensor decompositions by order comparison. Let next
$M=t^{1+2m}.\Sp_{2m}(t)\bZ(G)$ be the normaliser of an extra-special
$t$-subgroup where $n=m^t$. As for the case of $\SL_n(q)$, a maximal abelian
subgroup in $M$ has size at most $t^{m(m+1)/2+m+1}\gcd(n,q+1)$, while $N$
contains a maximal torus of order at least $(q-1)^{n-1}$, respectively
$3^{(n-1)/2}$ if $q=2$. Comparing the orders we arrive only at the case
$2^{1+8}.\Sp_8(2)<\SU_{16}(3)$. Here, the $\ell$-parts of $M$ and $G$ only
agree for $\ell=17$, but then the centralisers of $\ell$-elements in $G$ are
too large.

Similarly, by slight variations of the considerations for $\SL_n(q)$, none of
the other types of geometric maximal subgroups apart from $M=\Sp_{2n}(q)$ can
contain~$N$; for the latter
we arrive at the condition $n=4=d$ for $M$ to contain an $\ell$-subgroup of
sufficient rank, and furthermore $q\le3$ for it to contain the normaliser of a
Sylow $d$-torus.
\end{proof}

\begin{prop}   \label{prop:S in SU}
 In the situation of Proposition~\ref{prop:over SU} assume $M$ is a maximal
 subgroup of $G$ in class $\cS$ containing $N$. Then $M<G$ is one of:
 \begin{enumerate}[\rm(1)]
  \item $\PSL_2(7)<\SU_3(3)$, with $d=6$, $\ell=7$;
  \item $3.\fA_6.2_3<\SU_3(5)$, with $d=2$, $\ell=3$;
  \item $3.\fA_7<\SU_3(5)$, with $d=2$, $\ell=3$, or $d=6$, $\ell=7$;
  \item $4\circ2.\fA_7<\SU_4(3)$, with $d=4$, $\ell=5$, or $d=6$, $\ell=7$;
  \item $4_2.\PSL_3(4)<\SU_4(3)$, with $d=6$, $\ell=7$;
  \item $\PSL_2(11)<\SU_5(2)$, with $d=10$, $\ell=11$;
  \item $3.M_{22}<\SU_6(2)$, with $d=10$, $\ell=11$; or
  \item $3_1.\PSU_4(3).2_2<\SU_6(2)$, with $d=3$, $\ell=7$.
 \end{enumerate}
\end{prop}

\begin{proof}
First assume $S:=F^*(M)/\bZ(F^*(M))$ is not of Lie type in the same
characteristic as $G$. For $n\le12$ we again extract the relevant cases from
the tables in \cite[\S8]{BHR}. Note that $|\SU_n(q)|>|\SL_n(q)|$ for $n\ge3$.
Then arguing as in the proof of Lemma~\ref{lem:N in SL} we obtain that
$|N|\ge q^{n-1}$ here as well. The possibilities for $S$ not excluded by this
lower bound are now handled using the Atlas, leading to Cases~(1)--(8). We may
now assume $n\ge13$, and so $N$ contains a torus of size at least
$(q^2-1)^6\ge3^6$ if $q\le3$, and $3^{11}$ if $q\ge4$.

The considerations in the proof of Proposition~\ref{prop:S in SL} for $S$
alternating apply verbatim to show that no embeddings in dimension $n\ge13$
lead to examples.
In fact, the same is true for the whole discussion of cross characteristic
embeddings of groups of Lie type as well as for embeddings of sporadic groups,
again showing that no case with $n\ge13$ arises.

If $S$ is of Lie type in defining characteristic, we use again \cite[\S8]{BHR}
when $n\le12$ to rule out the occurrence of an example. For $n\ge13$ we can
argue as in the case of $\SL_n(q)$ that the $\ell$-rank of any candidate $S$
is too small.
\end{proof}

\begin{rem}   \label{rem:uniq SU}
 Again we see that in the situation of Proposition~\ref{prop:over SU} the
 normaliser of a Sylow $d$-torus lies in a \emph{unique} maximal subgroup of
 $\SU_n(q)$ whenever Sylow $\ell$-subgroups are non-cyclic, except in
 $\SU_3(5)$ with $\ell=3$ (where Sylow 3-subgroups are non-abelian).
 (Again, compare to \cite{BBGT}.)
\end{rem}

The $d$-split Levi subgroups of $G=\SU_n(q)$ have the form
$$\Big(\GL_{n_1}((-q)^e)\times\cdots\times\GL_{n_t}((-q)^e)\times\GU_s(q)\Big)
  \cap G\quad\text{with $e\sum_i n_i+s=n$}$$
by \cite[Exmp.~3.5.14]{GM20}, where we may and will assume
$n_1\ge\ldots\ge n_t$; note that here $s\equiv r\pmod e$ and thus $s\ge r$.
For abbreviation we will write $L_e(n_1,\ldots,n_t;s)$ for such a Levi
subgroup.

With this information in place we can determine the subnormalisers:

\begin{thm}   \label{thm:subn SU}
 Let $G=\SU_n(q)$ with $n\ge3$, let $\ell\nmid q$ be a prime such that Sylow
 $\ell$-subgroups of $G$ are abelian and $\bS_d\le\bG$ a Sylow $d$-torus,
 where $d:=e_\ell(q)$. With $e:=e_\ell(-q)$ write $n=ae+r$ with $0\le r<e$.
 Then for $x\in\bS_d^F$ an $\ell$-element we have one of
 \begin{enumerate}[\rm(1)]
  \item $\Sub_G(x)=\bN_G(\bS_d)$ if $\bC_G(x)=\bC_G(\bS_d)$;
  \item $\Sub_G(x)=(\GU_{ae}(q)\times\GU_r(q))\cap G$ if $r>0$ and
   $\bC_G(x)=L_e(n_1,\ldots,n_t;r)$ with $a=\sum n_i$ and at least one $n_i>1$;
  \item $\Sub_G(x)=(\GU_e(q)\wr\fS_a)\cap G$ if $r=0$, $e>1$, $a>1$ and
   $\bC_G(x)= L_e(1,\ldots,1;e)$; or
  \item $\Sub_G(x)=G$ otherwise.
 \end{enumerate}
\end{thm}

\begin{proof}
We proceed as in the proof of Theorem~\ref{thm:subn SL}. Observe that
$\ell>2$ by our assumptions. By Proposition~\ref{prop:cyclic} the cyclic Sylow
case again leads to~(1). Now, for all maximal subgroups in class $\cS$ coming
up in Proposition~\ref{prop:S in SU}, the Sylow $\ell$-subgroups of $G$ are
either cyclic or non-abelian, so no further cases arise from these.

Now assume $\Sub_G(x)$ lies in one of the maximal subgroups $M$ listed in
Proposition~\ref{prop:over SU}. Again by \cite[Prop.~2.2]{CE94}, $\bC_G(x)$ is
a $d$-split Levi subgroup of $G$ and hence has the form $L_e(n_1,\ldots,n_t;s)$
introduce above, where either $s>r$ or $n_1>1$ as otherwise
$\bC_G(x)=\bC_G(\bS_d)$ and we are in~(1) of the conclusion.
Now the groups in Cases~(3), (4) and~(5) in 
Proposition~\ref{prop:over SU} have cyclic Sylow $\ell$-subgroup. If we are in
Case~(2), and so $r=0$, then arguing as for $\SL_n(q)$ we see that
$\Sub_G(x)\le M$ as in~(3) of the conclusion if and only if $\bC_G(x)$ is as
claimed. Note here that when $e=1$ we obtain~(1), and when $a=1$ we have
$\Sub_G(x)=G$ as in~(4). Finally, if $M$ is a Levi subgroup of~$G$ as in
Case~(1) of Proposition~\ref{prop:over SU} then necessarily $r=s$ for
$\bC_G(x)$ and $\Sub_G(x)=\langle\bC_G(x),\bN_G(\bS_d)\rangle=M$ as in~(2) of
the conclusion by Proposition~\ref{prop:Burnside}. 
\end{proof}

\subsection{The symplectic and odd-dimensional orthogonal groups}
We now turn to $G$ either $\Sp_{2n}(q)$ or $\SO_{2n+1}(q)$ with $n\ge2$,
$(n,q)\ne(2,2)$. Let $d\ge1$ be an integer and
$$e:=\begin{cases} d& \text{if $d$ is odd},\\ d/2& \text{if $d$ is even}.
  \end{cases}$$
We write $n=ae+r$ with $0\le r<e$. Then the centraliser and normaliser of a
Sylow $d$-torus $\bS_d$ of the underlying simple algebraic group $\bG$ are
given as follows (see \cite[Exmp.~3.5.15 and 3.5.29]{GM20}):
$$\bC_G(\bS_d)=\GL_1(q^e)^a\times H_r\quad\text{and}\quad
  \bN_G(\bS_d)=\GL_1(q^e)\wr G(2e,1,a)\times H_r$$
if $d=e$ is odd, and
$$\bC_G(\bS_d)=\GU_1(q^e)^a\times H_r\quad\text{and}\quad
  \bN_G(\bS_d)=\GU_1(q^e)\wr G(2e,1,a)\times H_r$$
if $d=2e$ is even, where $H_r$ is a group of rank~$r$ of the same classical
type as $G$. We further assume there is a Zsigmondy prime divisor $\ell>2$ of
$q^d-1$, so that $d=e_\ell(q)$, $e=e_\ell(q^2)$.

\begin{prop}   \label{prop:over Sp}
 Let $G=\Sp_{2n}(q)$ with $n\ge2$ and $d,e,\ell$ as above. If $M<G$ is a
 maximal subgroup containing the normaliser of a Sylow $d$-torus of $\bG$,
 then one of:
 \begin{enumerate}[\rm(1)]
  \item $M=\Sp_{2ae}(q)\times\Sp_{2r}(q)$ if $r>0$;
  \item $M=\Sp_{2e}(q)\wr\fS_a$ if $n=ae$ and $a>1$;
  \item $M=\GL_n(q).2$ if $n=d$ is odd and $q$ is odd;
  \item $M=\Sp_{2n/t}(q^t).(q-1,2,t)$ if $n=e$ and $t|n$ is a prime;
  \item $M=\GU_n(q).2$ if $n=e=d/2$ is odd and $q$ is odd;
  \item $M=\Sp_{2n}(2)$ if $q=4$ and $d=1$;
  \item $M=\GO_{2n}^+(q)$ if $q$ is even, $n=ae>2$ and $2n/d$ is even;
  \item $M=\GO_{2n}^-(q)$ if $q$ is even, $n=ae>2$ and $2n/d$ is odd; or
  \item $M$ is in class $\cS$.
 \end{enumerate}
\end{prop}

\begin{proof}
Let $N:=\bN_G(\bS_d)$ for a Sylow $d$-torus $\bS_d\le\bG$.
By the description recalled above, $N$ acts irreducibly on subspaces of
dimensions $ae, ae$ and $2r$ of the natural module of~$G$. We refer to
\cite[Tab.~3.5.C]{KL90} for the Aschbacher classes.

The maximal parabolic subgroups $P_r$ (for $r>0$) do not contain a Sylow
$d$-torus of $G$ centralised by an $\Sp_{2r}(q)$-factor. Thus among reducible
maximal subgroups we only obtain the one in~(1). For the maximal imprimitive
subgroups $\Sp_{2n/t}(q)\wr\fS_t$ with $2\le t|n$, looking at the possible
imprimitivity decompositions for $N$ we arrive at the case in~(2). The
imprimitive subgroup $M=\GL_n(q).2$ with $q$ odd contains an elementary
abelian $\ell$-subgroup of the required rank only when $d$ is odd. Its
automiser is $G(2d,1,a)$ in $G$, and $G(d,1,a).2$ in $M$ by the description in
Section~\ref{subsec:SL}, which agree only when $a=1$. Comparing the
centralisers in $M$ and $G$ we see that then necessarily $d=n$, as in~(3) of
our conclusion.

The extension field subgroups $M=\Sp_{2m}(q^t).(q-1,2,t)$, with $n=mt$ and $t$
prime contain an elementary abelian $\ell$-subgroup $E$ of the right rank only
when $t|d$, respectively if $4|d$ when $t=2$. Now the automiser of $E$ in $M$
is $G(2e/t,1,a).(q-1,2,t)$ and $G(2e,1,a)$ in $G$, so we need $a=1$. Moreover,
comparing centralisers we see $r=0$, whence $n=e$ as in~(4). A subgroup
$M=\GU_n(q).2$, with $q$ odd, contains $E$ only when $d\equiv2\pmod4$, so
$e=d/2$ is odd, by the order formula. The automiser of $E$ in $M$ is now
$G(e,1,a).2$, and $G(d,1,a)$ in $G$, forcing $a=1$. Comparison of centralisers
shows $r=0$ and so $n=e$, as in~(5).
For the subfield subgroups $\Sp_{2n}(q_0)$ with $q=q_0^t$, $t$ prime, we can
argue as for $\SL_n(q)$ to arrive at $q=4$ and $d=1$, which in this case does
lead to a containment, listed in~(6).

The tensor product stabilisers of type $\Sp_{2m}(q)\GO_{n/m}^{(\pm)}(q)$, with
$m|n$, $n/m\ge3$ and $q$ odd, have strictly smaller $\ell$-rank than $G$ for
all relevant primes~$\ell$. For $M$ the normalizer of an extra-special subgroup
of type $2^{1+2m}.\GO_{2m}^-(2)$, with $n=2^{m-1}$ and $q=p$ odd, the size of
an abelian subgroup of $M$ containing our elementary abelian $\ell$-subgroup
is at most $3\cdot2^{m(m-1)/2+m+1}$ by \cite[Thm~3.1]{Vd99}, while it is at
least $(q-1)^n$ in $G$, respectively $7^{n/2}$ if $q=3$. This forces $m=4$ and
$q=3$, so $M=2^{1+8}.\GO_8^-(2)$ in $G=\Sp_{16}(3)$. Here, only $\ell=17$ is
possible, but elements of order~17 in~$G$ have too large a centraliser, so no
example arises here.

Comparing orders it is easily seen that the tensor induced subgroups
$\Sp_m(q)\wr\fS_t$ with $2n=m^t$, $t\ge3$ odd, cannot contain an elementary
abelian $\ell$-subgroup of the required rank. Finally, for $M=\GO_{2n}^\pm(q)$
with $q$ odd, comparing centralisers of elementary abelian subgroups in
Sylow $d$-tori, we arrive at the stated conditions in~(7) and~(8). In these
cases, $M$ does indeed contain a conjugate of~$N$ as can be seen from the
description of the Sylow $d$-normalisers in orthogonal groups recalled
in Section~\ref{subsec:GO} below.
\end{proof}

We need not discuss the orthogonal groups $\SO_{2n+1}(q)$ for even $q$, since
these are isomorphic to $\Sp_{2n}(q)$, nor $\SO_5(q)$ which possesses the same
non-abelian simple composition factor as $\Sp_4(q)$ and otherwise only
composition factors of order~2, whence subnormalisers of $\ell$-elements, for
$\ell\ne2$, of the two groups determine each other by \cite[Lemma~2.15]{Ma25}.

\begin{prop}   \label{prop:over SO}
 Let $G=\SO_{2n+1}(q)$ with $n\ge3$ and $q$ odd, and $d,e,\ell$ as above. If
 $M<G$ is a maximal subgroup containing the normaliser of a Sylow $d$-torus of
 $\bG$, then one of:
 \begin{enumerate}[\rm(1)]
  \item $M=(\GO_{2ae}^+(q)\times\GO_{2r+1}(q))\cap G$ if $2ae/d$ is even;
  \item $M=(\GO_{2ae}^-(q)\times\GO_{2r+1}(q))\cap G$ if $2ae/d$ is odd; or
  \item $M$ is in class $\cS$.
 \end{enumerate}
\end{prop}

\begin{proof}
Let $N:=\bN_G(\bS_d)$ for a Sylow $d$-torus $\bS_d\le\bG$. Here, by the
structure recalled above, $N$ has irreducible constituents of dimensions
$ae, ae$ and $2r+1$ on the natural module. We refer to
\cite[Tab.~3.5.D]{KL90} for the Aschbacher classes.

The maximal parabolic subgroups $P_r$ (for $r>0$) do not possess an
$\SO_{2r+1}(q)$-factor centralising a Sylow $d$-torus of $G$, so the only
reducible cases are those in (1) and~(2). By the order formulas, none of the
other types of geometric subgroups can contain our subgroup~$N$. For this,
observe that all but the subfield subgroups are of strictly smaller rank
than $G$.
\end{proof}

\begin{lem}   \label{lem:|Sp_n|}
 Let $n\ge2$ and $q$ be a prime power. Then
 $$|\Sp_{2n}(q)|=|\SO_{2n+1}(q)|>(q^2-1)^2q^{2n^2+n-4}.$$
\end{lem}

\begin{proof}
The order formula shows that $|\Sp_{2n}(q)|_{q'}=(q^2-1)|\SL_n(q^2)|_{q'}$,
thus an application of Lemma~\ref{lem:|SL_n|} gives the claim.
\end{proof}

\begin{lem}   \label{lem:N in Sp}
 For $G=\Sp_{2n}(q)$ or $\SO_{2n+1}(q)$ with $n\ge2$ the subgroup $N$ has
 order $|N|\ge 2^a q^n$, where $n=ae+r$ with $0\le r<e$.
\end{lem}

\begin{proof}
Arguing as in the proof of Lemma~\ref{lem:N in SL} we see
$$(q^e\pm1)^a|G(2e,1,a)|\ge ((3/4)q)^e(2e)^a a!\ge 2^aq^{ae}$$
and then Lemma~\ref{lem:|Sp_n|} yields the stated lower bound.
\end{proof}

We phrase the next result for the derived subgroup
$\Om_{2n+1}(q)=[\SO_{2n+1}(q),\SO_{2n+1}(q)]$ of index~2 in $\SO_{2n+1}(q)$
since there the normalisers of simple subgroups are more easily understood.
As all Sylow $d$-tori are conjugate,
$|\bN_{\SO_{2n+1}(q)}(\bS_d):\bN_{\Om_{2n+1}(q)}(\bS_d)|$ as well.

\begin{prop}   \label{prop:S in Sp}
 Let $G=\Sp_{2n}(q)$ with $n\ge2$, or $G=\Om_{2n+1}(q)$ with $n\ge3$ and $q$
 odd and assume $M\in\cS$ contains the normaliser $N$ of a Sylow $d$-torus.
 Then $M<G$ is one of:
 \begin{enumerate}[\rm(1)]
  \item $\PSU_3(3).2<\Sp_6(2)$, with $d=3$, $\ell=7$;
  \item $2.\PSL_2(13)<\Sp_6(3)$, with $d=3$, $\ell=13$;
  \item $\PSL_2(17)<\Sp_8(2)$, with $d=8$, $\ell=17$;
  \item $\fS_{10}<\Sp_8(2)$, with $d=4$, $\ell=5$, or $d=3$, $\ell=7$;
  \item $\fS_{14}<\Sp_{12}(2)$, with $d=3$, $\ell=7$;
  \item $\fS_9<\Om_7(3)$, with $d=4$, $\ell=5$, or $d=6$, $\ell=7$; or
  \item $G_2(3)<\Om_7(3)$, with $d=3$, $\ell=13$.
 \end{enumerate}
\end{prop}

\begin{proof}
Let $S:=F^*(M)/\bZ(F^*(M))$. Again we start with the case that $S$ is not of
Lie type in characteristic dividing~$q$. As before, the case $n\le6$ for the
symplectic groups respectively $n\le5$ for the orthogonal groups can be
discussed using the tables in \cite[\S8]{BHR}, which leads to the groups
in~(1)--(6) of our conclusion. So we will assume now that $n\ge7$
(resp.~$n\ge6$) and thus $N$ contains a maximal (abelian) torus of order at
least $(q-1)^6\ge 3^6$ when $q\ge4$, at least $(q^2-q+1)^3=7^3$ when
$q=3$ and at least $(q^2-1)^4=3^4$ when $q=2$ (using that $q$ is odd when
$G=\Om_{2n+1}(q)$).

First assume $S=\fA_m$. Since $2n+1\ge13$ we then have $m\le 2n+2$. The
normaliser $N$ contains a maximal torus $T$ of order at least $(q-1)^n$, so at
least $3^n$ if $q\ge4$. Since maximal abelian subgroups of $\fS_m$ have size at
most $3^{m/3}\le 3^{(2n+2)/3}$, comparison shows we must have $q\le3$.
If $q=3$ then $|T|=(q^e\pm1)^a(q+1)^r\ge 8^{n/2}$ and comparing with the bound
for $\fS_m$ we arrive at $n\le 2.4<6$. Assume $q=2$, and so $G$ is symplectic.
If $\ell=3$, so $d=2$,
then $|T|=3^n$ which is always larger than $3^{(2n+2)/3}$. When $\ell\ge5$
observe that $\ell$-elements in $T$ are conjugate to $2e$ of their powers,
while $\ell$-elements in $\fS_m$ are rational, so we conclude that $\ell=2e+1$. 
Now a maximal elementary abelian $\ell$-subgroup $E$ of $\fS_m$ has rank
$\lfloor m/\ell\rfloor$ and the largest abelian subgroup containing it has
order at most $3^{m'}|E|$ with $m'=(m-\ell\lfloor m/\ell\rfloor)/3$, while in
$G$ such a subgroup has rank $a=\lfloor n/e\rfloor$ and is contained in a torus
of order at least $(2^e-1)^a(2^r+1)$. Comparing the two we conclude that
$\ell=5$ and $7\le n\le 9$. In all of these cases, the centraliser of $E$ in
$\fS_m$ does not contain a subgroup $\Sp_{2r}(2)$.

If $S$ is sporadic, since $N$ contains a cyclic subgroup of order $q^e\pm1$,
as in the proof of Proposition~\ref{prop:S in SL} we conclude that $e\le6$,
and in fact $e\le4$ when $q\ge3$, and $e\le3$ when $q\ge4$. Since no covering
group of a sporadic simple group possesses elements of order~65, 80 or~82, we
have $e\le5$ when $q=2$ and
$e\le3$ when $q\ge4$. Now when $e\le3$ then $a\ge2$ as $n\ge7$ and hence
$\ell^2$ divides $|S|$. In this case, as well as when $q=2$, $e=4,5$ and
$7\le n\le9$, the table of low dimensional representations in \cite{HM} shows
that no example arises.

Next assume $S$ if of Lie type in cross characteristic. For $S=\PSL_m(y)$ with
$m\ge3$, as in the proof of Proposition~\ref{prop:S in SL} the smallest degree
is $(y^m-1)/(y-1)-m\ge y^{m-1}$, unless
$(m,y)\in \cE:=\{(3,2),(3,4),(4,2),(4,3)\}$, and we have
$|M|\le|\Aut(S)|\le 2|\PGL_m(y)|\log y\le 2\ y^{m^2-1}\log y$.
Thus, for $(m,y)\notin \cE$, if $N\le M$ then
$$2q^{y^{m-1}/2}\le 2^a q^n\le |N|\le |M|\le 2 y^{m^2-1}\log y$$
by Lemma~\ref{lem:N in Sp}. Assuming $2n+1\ge13$, this is only satisfied for
$S=\PSL_3(5),\PSL_5(2),\PSL_6(2)$, with $2n+1$ at most 19, 14, 21.
But none of these groups has a faithful projective representation in this
range, by \cite{HM}. The groups with $(m,y)\in \cE$ can be excluded in a
similar way, so no further case arises.

For $S=\PSL_2(y)$, $y\ge7$ and $y\ne9$, the Sylow $\ell$-subgroups of $S$ are
cyclic, with centraliser order at most $y+1$, while the centraliser of an
$\ell$-element in $T\le G$ has size at least $q^{(n+1)/2}-1$, so
$y\ge q^{(n+1)/2}-2$. On the other hand, the minimal projective degree of $S$
is $(y-1)/(2,y-1)$, so we need $n\ge (y-1)/2$.
This is only satisfied for values $n\le7$ which does not lead to new cases.
The other series of groups can be dealt with by analogues estimates.

If $S$ is of Lie type in characteristic $p$, the low dimensional cases can
again be discussed using \cite[\S8]{BHR}, which leads only to the example
in~(7) of the conclusion. So assume $n\ge7$, respectively $n\ge6$ for
$G=\Om_{2n+1}(q)$. If $S=\PSL_m(p^f)$ then by \cite[Prop.~5.4.6, 5.4.8]{KL90}
we must have $2n+1\ge m(m-1)$ (since the representation needs to be self-dual)
and $p^f|q$, or $2n+1\ge m^k$ and
$p^f|q^k$ for some $k\ge2$. Write $m=bd+s$, so $S$ has $\ell$-rank at
most $b$, then $2n+1\ge (bd+s)(bd+s-1)$ and since $d\ge e$ we find
$n\ge e(eb^2/2+bs-b/2)$ whence $G$ has $\ell$-rank at least $eb^2/2+bs-b/2$.
This is larger than $b$ unless $m\le6$, and in those cases we obtain $n\le5$,
which was excluded. Again the same type of estimates applies to all other types
of $S$, leading to no further candidates.
\end{proof}

\begin{rem}   \label{rem:uniq Sp}
(1) The maximal subgroups of $\Om_7(3)$ in (6) and (7) of the preceding
 result are not stable under the diagonal automorphism induced by
 $\SO_7(3)$ (see \cite[Tab.~8.40]{BHR}).   \par
(2) The above shows that for the symplectic groups $\Sp_{2n}(q)$ with $q$ even
 it may happen that the normaliser of a Sylow $d$-torus lies in two distinct
 maximal subgroups even when Sylow $\ell$-subgroups are non-cyclic, namely in
 Cases~(2) and~(7), (8) of Proposition~\ref{prop:over Sp}. Further such cases
 occur for $\ell\le7$. For $\SO_{2n+1}(q)$ with $q$ odd, there is always a
 unique maximal overgroup in the non-cyclic Sylow case.
\end{rem}

We next determine the subnormalisers. For this we recall from
\cite[Exmp.~3.5.15]{GM20} that the $d$-split Levi subgroups of $G$ have the
form
$$\GL_{n_1}(q^e)\times\cdots\times\GL_{n_t}(q^e)\times H_s\qquad
  \text{if $d=e$ is odd},$$
respectively
$$\GU_{n_1}(q^e)\times\cdots\times\GU_{n_t}(q^e)\times H_s\qquad
  \text{if $d=2e$ is even},\quad$$
where $H_s=\Sp_{2s}(q)$ resp.\ $\SO_{2s+1}(q)$ and
$e\sum_i n_i+s=n$ in either case, where we may assume
$n_1\ge n_2\ge\cdots\ge n_t$. Here again $s\equiv r\pmod e$ and thus $s\ge r$.
We will write $L_d(n_1,\ldots,n_t;s)$ for such a Levi subgroup (whose structure
depends on the parity of $d$ as indicated).

\begin{thm}   \label{thm:subn Sp}
 Let $G=\Sp_{2n}(q)$ with $n\ge2$, let $\ell\nmid q$ be a prime such that Sylow
 $\ell$-subgroups of $G$ are abelian and $\bS_d\le\bG$ a Sylow $d$-torus,
 where $d:=e_\ell(q)$. With $e:=e_\ell(q^2)$ write $n=ae+r$ where $0\le r<e$.
 Then for $x\in\bS_d^F$ an $\ell$-element we have one of
 \begin{enumerate}[\rm(1)]
  \item $\Sub_G(x)=\bN_G(\bS_d)$ if $\bC_G(x)=\bC_G(\bS_d)$;
  \item $\Sub_G(x)=\Sp_{2ae}(q)\times\Sp_{2r}(q)$ if $r>0$, $q$ is odd and
   $\bC_G(x)=L_d(n_1,\ldots,n_t;r)$ with $n_1>1$;
  \item $\Sub_G(x)=\GO_{2ae}^+(q)\times\Sp_{2r}(q)$ if $q$ is even,
   $2ae/d$ is even and $\bC_G(x)=L_d(n_1,\ldots,n_t;r)$ with $n_1>1$;
  \item $\Sub_G(x)=\GO_{2ae}^-(q)\times\Sp_{2r}(q)$ if $q$ is even,
   $2ae/d$ is odd and $\bC_G(x)=L_d(n_1,\ldots,n_t;r)$ with $n_1>1$;
  \item $\Sub_G(x)=\Sp_{2e}(q)\wr\fS_a$ if $r=0$, $a>1$ and
   $\bC_G(x)= L_d(1,\ldots,1;e)$; or
  \item $\Sub_G(x)=G$ otherwise.
 \end{enumerate}
\end{thm}

\begin{proof}
The cyclic Sylow case is again covered by Proposition~\ref{prop:cyclic}.
So assume Sylow $\ell$-subgroups of $G$ are non-cyclic and $\Sub_G(x)$ is
contained in a maximal subgroup $M$ of $G$. Among the maximal subgroups in
Proposition~\ref{prop:S in Sp} only $\fS_{10}<\Sp_8(2)$ with $\ell=5$ and
$\fS_{14}<\Sp_{12}(2)$ with $\ell=7$ possess 
non-cyclic Sylow $\ell$-subgroups. Direct computation in \GAP\ shows that in
either case, the subnormalisers of $\ell$-elements are as in~(3) or~(5) of
the conclusion (and hence not symmetric groups).

Now assume $\bC_G(x)$, and hence $\Sub_G(x)$, is contained in one of the
maximal subgroups $M$ in Proposition~\ref{prop:over Sp}. Again by
\cite[Prop.~2.2]{CE94} the centraliser $\bC_\bG(x)$ is a $d$-split Levi
subgroup of $\bG$. If $r=0$ then Cases~(2)--(8) in
Proposition~\ref{prop:over Sp} are relevant. In Cases~(3), (4) and~(5) the
Sylow $\ell$-subgroups of $G$ are cyclic. In Case~(6) we have $\ell=3$ and
Sylow $3$-subgroups are abelian only when $n=2$, a possibility that can be
ruled out by explicit computation. If $\bC_G(x)$ lies in
$M=\Sp_{2e}(q)\wr\fS_a$ then as in the linear and unitary cases necessarily
$s=e$, while all $n_i=1$, giving~(5). If $\bC_G(x)$ lies in an orthogonal group
as in Cases~(7) or~(8) then comparing centralisers shows that $s=0$, leading to
the conclusion in (3) or~(4).

Finally, if $r>0$ then the only maximal overgroup for $N$ is a subsystem
subgroup $M=\Sp_{2ae}(q)\times\Sp_{2r}(q)$ as in Case~(1) of
Proposition~\ref{prop:over Sp}. Comparing dimensions of composition factors
on the natural module we conclude that $r=s$, and in this case, since $M$
contains the normaliser of a Sylow $d$-torus, $\Sub_G(x)\le M$ by
Proposition~\ref{prop:Burnside}. Now going again through the cases in
Proposition~\ref{prop:over Sp} for the $\Sp_{2ae}(q)$-factor it follows that
we must in fact have $\Sub_G(x)=M$ when $q$ is odd, while when $q$ is even,
$\bC_G(x)$ and hence $\Sub_G(x)$ lies in a subgroup as given in~(3) or~(4).
\end{proof}

\begin{thm}   \label{thm:subn SOodd}
 Let $G=\SO_{2n+1}(q)$ with $n\ge3$ and $q$ odd, let $\ell\nmid q$ be a prime
 such that Sylow $\ell$-subgroups of $G$ are abelian and $\bS_d\le\bG$ a Sylow
 $d$-torus, where $d:=e_\ell(q)$. With $e:=e_\ell(q^2)$ write $n=ae+r$ where
 $0\le r<e$. Then for $x\in\bS_d^F$ an $\ell$-element one of
 \begin{enumerate}[\rm(1)]
  \item $\Sub_G(x)=\bN_G(\bS_d)$ if $\bC_G(x)=\bC_G(\bS_d)$;
  \item $\Sub_G(x)=\big(\GO_{2ae}^+(q)\times\GO_{2r+1}(q)\big)\cap G$ if
   $2ae/d$ is even and $\bC_{G}(x)=L_d(n_1,\ldots,n_t;r)$ with $n_1>1$;
  \item $\Sub_G(x)=\big(\GO_{2ae}^-(q)\times\GO_{2r+1}(q)\big)\cap G$ if
   $2ae/d$ is odd and $\bC_{G}(x)=L_d(n_1,\ldots,n_t;r)$ with $n_1>1$; or
  \item $\Sub_G(x)=G$ otherwise.
 \end{enumerate}
\end{thm}

\begin{proof}
The proof is very similar to the one of Proposition~\ref{thm:subn Sp}, but
easier, as there are fewer cases to consider: The two groups in
Proposition~\ref{prop:S in Sp}(6) and~(7) possess cyclic Sylow $\ell$-subgroups
and thus do not occur as subnormalisers by Proposition~\ref{prop:cyclic}
while the maximal subgroups in~(1) and~(2) of Proposition~\ref{prop:over SO}
lead to~(2) and~(3) of the conclusion.
\end{proof}

\subsection{The even-dimensional orthogonal groups}   \label{subsec:GO}
Finally, we consider the even-dimen\-sional special orthogonal groups
$\SO_{2n}^\pm(q)$, $n\ge4$. Here, by convention we let
$\SO_{2n}:=\GO_{2n}^\circ$, the connected component of the identity, and
$\SO_{2n}^\pm(q)$ the group of fixed points under a Frobenius map $F$ with
respect to an $\FF_q$-structure.

Let $d\ge1$ be an integer and as before set $e:=d$ if $d$ is odd, $e:=d/2$ if
$d$ is even. We write $ n=ae+r$ with $0\le r<e$, except if $d|2n$ and either
$2n/d$ is odd in type $\SO_{2n}^+(q)$, or $2n/d$ is even in type
$\SO_{2n}^-(q)$: in the latter cases, write $n=(a+1)e$ and set $r:=e$ (note
that $e|n$ under our assumptions). Then the centraliser and normaliser of a
Sylow $d$-torus $\bS_d$ of~$\bG=\SO_{2n}$ in $\hat G:=\GO_{2n}^\pm(q)$ are
given as follows (see \cite[Exmp.~3.5.15 and 3.5.29]{GM20}):
$$\bC_{\hat G}(\bS_d)=\GL_1(q^e)^a\times H_r\quad\text{and}\quad
  \bN_{\hat G}(\bS_d)=\GL_1(q^e)\wr G(2e,1,a)\times H_r$$
if $d=e$ is odd, and
$$\bC_{\hat G}(\bS_d)=\GU_1(q^e)^a\times H_r\quad\text{and}\quad
  \bN_{\hat G}(\bS_d)=\GU_1(q^e)\wr G(2e,1,a)\times H_r$$
if $d=2e$ is even, where $H_r$ is a general orthogonal group of
rank~$r$ of the same type as $G$, except that it is of opposite sign if $d$ is
even and $a$ is odd. Again, we assume there is a Zsigmondy prime divisor
$\ell>2$ of $q^d-1$, so that $d=e_\ell(q)$, $e=e_\ell(q^2)$.

We again first consider the maximal overgroups of Sylow tori normalisers:

\begin{prop}   \label{prop:over SO+}
 Let $G=\SO_{2n}^+(q)$ with $n\ge4$, and $d,e,\ell$ as above. If $M<G$ is a
 maximal subgroup containing the normaliser of a Sylow $d$-torus,
 then one of:
 \begin{enumerate}[\rm(1)]
  \item $M=(\GO_{2ae}^+(q)\times\GO_{2r}^+(q))\cap G$ if $r>0$, and $a$ is
   even or $d$ is odd;
  \item $M=(\GO_{2ae}^-(q)\times\GO_{2r}^-(q))\cap G$ if $r>0$, $a$ is odd and
   $d$ is even;
  \item $M=\Sp_{2n-2}(2)$ if $q=2$, $r=1$, and $a$ is even or $d>1$ is odd;
  \item $M=(\GO_{2e}^+(q)\wr\fS_a)\cap G$ if $r=0$, $a>1$ and $d$ is odd;
  \item $M=(\GO_{2e}^-(q)\wr\fS_a)\cap G$ if $r=0$, $a$ is even and $d$ is even;
  \item $M=\GL_n(q).2$ if $e=d$ is odd, and $n=e+1$ or $n=2e$;
  \item $M=\GU_n(q).2$ if $e=d/2$ is odd, and $n=e+1$ or $n=2e$;
  \item $M=\GO_{2n/t}^+(q^t).t\cap G$ if $n=d$ is odd, $2<t|n$ is prime and
   $2n/t\ge4$;
  \item $M=\GO_n^+(q^2).2\cap G$ if $n=d\equiv0\pmod4$;
  \item $M=\GO_{2n}^+(2)$ if $q=4$, $d=1$ and $n$ is even;
  \item $M=\GO_{2n}^-(2)$ if $q=4$, $d=1$ and $n$ is odd; or
  \item $M$ is in class $\cS$.
 \end{enumerate}
\end{prop}

Observe that in items (4) and (5) for $e=1$, so $d\in\{1,2\}$, we have $M=N$,
hence $N$ is maximal in those cases.

\begin{proof}
Let $\bS_d$ be a Sylow $d$-torus of $(\bG,F)$. The structure of $\bC_G(\bS_d)$
and $N:=\bN_G(\bS_d)$ was recalled above. We refer to \cite[Tab.~3.5.E]{KL90}
and the relevant tables in \cite{BHR} for a description of the Aschbacher
classes. The composition factors of $N$ on the natural module have dimensions
$ae, ae$ and~$2r$. From this it follows that reducible maximal subgroups
containing $N$ are either as in~(1), (2), or, if $N$ lies in a subgroup
$M=\Sp_{2n-2}(q)$ with $q$ even, then we must have $2n-2=2ae$ and thus $r=1$.
In this case, the centraliser of $\bS_d^F$ in $G$ is larger than the one in $M$
by a factor $\SO_{2}^\eps(q)$, which forces $q=2$ and $\eps=+$, and the latter
holds when $a$ is even or $d$ is odd, giving~(3).

The maximal imprimitive subgroups lead to cases~(4)--(6). Here note that
normalisers of subgroups of type $\GL_n(q)$ are only maximal in $\SO_{2n}^+(q)$
when $n$ is even. In this case, they contain a Sylow $d$-torus of $G$ when
$d\ge3$ is odd, and its full normaliser only when moreover either $n=2d$ or
$n=d+1$. Further, the imprimitive maximal subgroup $\GO_n(q)\wr2$, with $nq$
odd, is of smaller rank than~$G$ and thus the centraliser of its Sylow
$d$-torus is smaller than the one in $G$ (using the precise descriptions in
\cite[Prop.~4.2.14, 4.2.16]{KL90}). The same holds for the extension field
subgroup $\GO_n(q^2)$, $nq$ odd, using \cite[Prop.~4.3.20]{KL90}.

Normalisers of subgroups of type $\GU_n(q)$, with $n$ even, contain a Sylow
$d$-torus only when $d\equiv2\pmod4$, and the full normaliser of $\bS_d$ only
when either $n=d$ or $n=d/2+1$, giving~(7).
For $M$ of type $\GO_{2n/t}^+(q^t)$, with $t|n$ a prime and $n/t\ge2$, to
contain a Sylow $d$-torus, we need $t|d$, and $4|d$ if $t=2$. Comparing
centralisers yields that $d$ must divide $n$. Now if $t$ is odd, the automiser
of $\bS_d$ in $M$ is $G(2e/t,2,a).t$ and $G(2e,2,a)$ in $G$, forcing $a=1$, so
$n=e$ must be odd (as $d|n$). If $t=2$ then $4|d|n$ forces $a\ge2$.
The automiser of $\bS_d$ in $M$ is $G(d/2,1,a).2$ but $G(d,2,a)$ in $G$, which
agree only when $a=2$, so $n=2e=d$. We hence reach the conditions stated in~(8)
and~(9).

For the subfield subgroups, comparing orders of maximal tori as in the
earlier proofs, we see that necessarily $q=4$ and $d=1$, and thus only the
groups in~(10) or~(11) can contain the normaliser of a Sylow $d$-torus. If $M$
is the normaliser of an extra-special 2-group, with $q$ odd, estimates as in
the proof for $\Sp_{2n}(q)$ rule out all cases with $m\ge4$ except
$M=2^{1+2m}.\GO_{2m}^+(2)<\SO_{2^m}^+(3)$ for $m=4,5$. For those, the only
prime $\ell$ such that $M$ contains a Sylow $\ell$-subgroup of $G$ is $\ell=31$
when $m=5$, but then its centraliser in $G$ is much too big.

Finally, as in the other types, the tensor product subgroups and the tensor
induced subgroups have too small a rank to contain an $\ell$-group of the
necessary rank.
\end{proof}

\begin{prop}   \label{prop:over SO-}
 Let $G=\SO_{2n}^-(q)$ with $n\ge4$ and let $d,e,\ell$ be as above. If $M<G$ is
 a maximal subgroup containing the normaliser of a Sylow $d$-torus, then one of:
 \begin{enumerate}[\rm(1)]
  \item $M=(\GO_{2ae}^+(q)\times\GO_{2r}^-(q))\cap G$ if $r>0$, and $a$ is even
   or $d$ is odd;
  \item $M=(\GO_{2ae}^-(q)\times\GO_{2r}^+(q))\cap G$ if $r>0$, $a$ is odd and
   $d$ is even;
  \item $M=\Sp_{2n-2}(2)$ if $q=2$, $r=1$, $a$ is odd and $d$ is even;
  \item $M=(\GO_{2e}^-(q)\wr\fS_a)\cap G$ if $r=0$, $a$ is odd and $d$ is even;
  \item $M=\GU_n(q)$ if $n=e=d/2$ is odd;
  \item $M=\GO_n(9)$ if $q=3$, $d\equiv0\pmod4$ and $n=d/2+1$;
  \item $M=\GO_{2n/t}^-(q^t).t\cap G$ if $n=e=d/2$, $t|n$ is prime and
   $2n/t\ge4$ is even; or
  \item $M$ is in class $\cS$.
 \end{enumerate}
\end{prop}

\begin{proof}
We use the description in \cite[Tab.~3.5.F]{KL90} of the Aschbacher classes.
Among reducible maximal subgroups we only find the examples listed in~(1)--(3),
with arguments entirely similar to the ones used for the groups of plus-type.
Next, using the dimensions of the $N$-composition factors on the natural
module, we find the imprimitive subgroups in~(4). Here, the
subgroup $\GO_n(q)^2$, with $nq$ odd, is of smaller rank than $G$ and using
the description in \cite[Prop.~4.1.6,4.2.16]{KL90} cannot contain a Sylow
$d$-torus normaliser.

The subgroup $M=\GU_n(q)$ with $n$ odd contains a suitable elementary abelian
$\ell$-subgroup only when $d\equiv2\pmod4$ and thus $e=d/2$ is odd. Its
automiser in $M$ is then $G(e,1,a)$, and in $G$ it equals $G(2e,2,a)$, so
we must have $a=1$. Comparing centralisers shows $r=0$, hence $n=e=d/2$ as
listed in (5). The extension field subgroups $\GO_n(q^2)$ with $nq$ odd do
only occur if $q=3$ and $d=2n-2$, as listed in~(6), taking into account the
information in \cite[Prop.~4.3.20]{KL90}, while the extension field subgroups
$\GO_{2n/t}^-(q^t).t$ with $t|2n$ prime, $2n/t\ge3$, can be discussed as in
the earlier proofs, leading to (7) of the conclusion. The tensor product
subgroups have smaller rank than necessary, and the maximal tori in subfield
subgroups of type $\GO_{2n}^-(q_0)$, with $q=q_0^t$ and $t\ge3$ prime, are too
small.
\end{proof}

In order to deal with the maximal subgroups in class $\cS$ we again first
derive a lower bound for $|N|$:

\begin{lem}   \label{lem:|SO_2n|}
 Let $n\ge4$ and $q$ be a prime power. Then
 $$|\SO_{2n}^\pm(q)|>\frac{1}{2}q^{2n^2-n}.$$
\end{lem}

\begin{proof}
Since $|\SO_{2n}^-(q)|>|\SO_{2n}^+(q)|$ it suffices to consider the latter
group. By the order formula $|\SO_{2n}^+(q)|=q^{n-1}(q^n-1)|\Sp_{2n-2}(q)|$
so from Lemma~\ref{lem:|Sp_n|} we get
$|\SO_{2n}^+(q)|>(q^2-1)^2(q^n-1)q^{2n^2-2n-4}$. Multiplying out we see
that $(q^2-1)^2(q^n-1)\ge (q-1)q^{n+3}$. Using that $q-1\ge q/2$ then achieves
the proof.
\end{proof}

\begin{lem}   \label{lem:N in SO}
 For $G=\SO_{2n}^\pm(q)$ with $n\ge4$ the subgroup $N$ has order
 $|N|\ge 2^{a-2} q^n$, where $n=ae+r$ with $0\le r\le e$ as introduced above.
\end{lem}

\begin{proof}
The relative Weyl group of a Sylow $d$-torus of $G$ is either $G(2e,1,a)$ or
its normal subgroup $G(2e,2,a)$ of index~2. Now arguing as in the proof of
Lemma~\ref{lem:N in SL} we see
$$(q^e\pm1)^a|G(2e,2,a)|\ge ((3/4)q)^e(2e)^a a!/2\ge 2^{a-1}q^{ae}.$$
Combining this with the bound in Lemma~\ref{lem:|SO_2n|} for $|\SO_{2r}^\pm(q)|$
we conclude.
\end{proof}

\begin{prop}   \label{prop:S in SO}
 Let $G=\Om_{2n}^\pm(q)$ with $n\ge4$, and assume $M$ is a maximal subgroup
 of~$G$ in class $\cS$ containing the normaliser of a Sylow $d$-torus. Then
 $M<G$ is one of:
 \begin{enumerate}[\rm(1)]
  \item $\fA_9<\Om_8^+(2)$, with $d=3$, $\ell=7$;
  \item $2.\Om_8^+(2)<\Om_8^+(3)$, with $d=4$, $\ell=5$;
  \item $\fA_{16}<\Om_{14}^+(2)$, with $d=3$, $\ell=7$; or
  \item $\fA_{12}<\Om_{10}^-(2)$, with $d=3$, $\ell=7$.
 \end{enumerate}
\end{prop}

\begin{proof}
Let $S=F^*(M)/\bZ(F^*(M))$. As before, the case $n\le6$ can be discussed using
the tables in \cite[\S8]{BHR}, which leads to Conclusions~(1), (2) and~(4),
so we may assume $n\ge7$. The sporadic groups are dealt with as in the proofs of
Propositions~\ref{prop:S in SL} and~\ref{prop:S in Sp} using \cite[Tab.]{HM}.
Next assume $S=\fA_m$ is alternating. As earlier this implies $m\le2n+2$.
Arguing exactly as in Proposition~\ref{prop:S in Sp} we see that necessarily
$q=2$ for both the plus- and the minus-type. Then again comparing the rank of a
maximal elementary abelian $\ell$-subgroup $E$ of $S$ and the size of a maximal
abelian subgroup containing it with the corresponding data in $G$, we find that
$\ell\le11$ and $n\le13$. Going through the cases, and using that
$$\fA_{2n+2}\le\begin{cases} \Om_{2n}^+(2)& \text{if $n\equiv3\pmod4$},\\
    \Om_{2n}^-(2)& \text{if $n\equiv1\pmod4$},\end{cases}$$
and
$$\fA_{2n+1}\le\begin{cases} \Om_{2n}^+(2)& \text{if $n\equiv0,3\pmod4$},\\
    \Om_{2n}^-(2)& \text{if $n\equiv1,2\pmod4$}\end{cases}$$
(see \cite[p.~187]{KL90}) we only arrive at the additional item~(3) of the
conclusion.

If $S$ is of Lie type in cross characteristic, we again use the lower bounds
for faithful projective representations from \cite{TZ96} and argue as before.
No further cases arise. The same holds for groups of Lie type in the same
characteristic.
\end{proof}

\begin{rem}
Observe that the maximal subgroups of type $\Om_8^+(2)$ of $\Om_8^+(3)$ are not
invariant under the outer diagonal automorphism induced by $\SO_8^+(3)$ (see
\cite[Tbl.~8.50]{BHR}) while a Sylow 4-torus is. Since all semisimple classes
are invariant under diagonal automorphisms, this means that $\Om_8^+(2)$ cannot
occur as a subnormaliser of a 5-element of $\SO_8^+(3)$.
\end{rem}

We now determine the subnormalisers. The $d$-split Levi subgroups of
$\GO_{2n}^\eps(q)$ have the form
$$\GL_{n_1}(q^e)\times\cdots\times\GL_{n_t}(q^e)\times\GO_{2s}^\eps(q)\qquad
  \text{if $d=e$ is odd},$$
respectively
$$\GU_{n_1}(q^e)\times\cdots\times\GU_{n_t}(q^e)\times\GO_{2s}^{\eps'}(q)\qquad
  \text{if $d=2e$ is even},\quad$$
where $\eps'=\eps(-1)^{\sum n_i}$ and $e\sum_i n_i+s=n$ in either case (see
\cite[Exmp.~3.5.15]{GM20}) and we may assume $n_1\ge\cdots\ge n_t$. Here again
we have $s\equiv r\pmod e$ and thus $s\ge r$. As before will write
$L_d(n_1,\ldots,n_t;s)$ for such a $d$-split Levi subgroup.

\begin{thm}   \label{thm:subn SOeven+}
 Let $G=\SO_{2n}^+(q)$ with $n\ge4$, let $\ell\nmid q$ be a prime such that
 Sylow $\ell$-subgroups of $G$ are abelian and $\bS_d\le\bG$ a Sylow $d$-torus
 where $d=e_\ell(q)$. Set $e:=e_\ell(q^2)$ and let $a,r$ be as defined above.
 Then for $x\in\bS_d^F$ an $\ell$-element we have one of
 \begin{enumerate}[\rm(1)]
  \item $\Sub_G(x)=\bN_G(\bS_d)$ if $\bC_G(x)=\bC_G(\bS_d)$;
  \item $\Sub_G(x)=\big(\GO_{2ae}^+(q)\times\GO_{2r}^+(q)\big)\cap G$ if $r>0$,
   $2ae/d$ is even, $\bC_{G}(x)=L_d(n_1,\ldots,n_t;r)$ with $n_1>1$;
  \item $\Sub_G(x)=\big(\GO_{2ae}^-(q)\times\GO_{2r}^-(q)\big)\cap G$ if $r>0$,
   $2ae/d$ is odd, $\bC_{G}(x)=L_d(n_1,\ldots,n_t;r)$ with $n_1>1$;
  \item $M=(\GO_{2e}^+(q)\wr\fS_a)\cap G$ if $r=0$, $a>1$, $d$ is odd and
   $\bC_{G}(x)=L_d(1,\ldots,1;e)$;
  \item $M=(\GO_{2e}^-(q)\wr\fS_a)\cap G$ if $r=0$, $a$ is even, $d$ is even
   and $\bC_{G}(x)=L_d(1,\ldots,1;e)$;
  \item $M=\GL_n(q).2$ if $e=d=n/2$ is odd and $\bC_{G}(x)=L_d(2;0)$;
  \item $M=\GU_n(q).2$ if $e=d/2=n/2$ is odd and $\bC_{G}(x)=L_d(2;0)$;
  \item $M=\GO_n^+(q^2).2\cap G$ if $n=d\equiv0\pmod4$ and
   $\bC_{G}(x)=L_d(2;0)$; or
  \item $\Sub_G(x)=G$ otherwise.
 \end{enumerate}
\end{thm}

\begin{proof}
With Proposition~\ref{prop:cyclic} we may assume Sylow $\ell$-subgroups are
non-cyclic and so $a>1$. The only groups in Proposition~\ref{prop:S in SO}
with non-cyclic Sylow $\ell$-subgroups are $2.\Om_8^+(2)$ in $\Om_8^+(3)$ with
$\ell=5$, and $\fA_{16}$ in $\Om_{14}^+(2)$ with $\ell=7$, but explicit
computation shows that in either case the subnormalisers are as in~(5) or~(2)
of the conclusion.

By \cite[Prop.~2.2]{CE94} the centraliser $C:=\bC_G(x)$ is $d$-split and thus
$C=L_d(n_1,\ldots,n_t;s)$ as introduced above. We may assume $C$ is contained
in one of the maximal subgroups $M$ of $G$ listed in
Proposition~\ref{prop:over SO+}
(as otherwise we have $\Sub_G(x)=G$ by Proposition~\ref{prop:Burnside}).
Also, if $C=\bC_G(\bS_d)$ then we are in Case~(1), so now suppose
$\bC_G(\bS_d)$ is strictly contained in $C$.
Assume first $r=0$. Then Cases~(4)--(9) are relevant; note that in Cases~(10)
and~(11) we have $\ell=3$ and Sylow $3$-subgroups of $G$ are non-abelian.
Also, in Case~(8) Sylow $\ell$-subgroups are cyclic, and similarly in Cases~(6)
and~(7) unless $n=2e$.
If $C$ lies in $\GO_{2e}^\pm(q)\wr\fS_a$ then necessarily $s=e$ and all
$n_i=1$, as in~(4) and~(5) of the conclusion. If $C$ lies in $\GL_n(\pm q).2$
then clearly $s=0$ and the assumptions in Proposition~\ref{prop:over SO+}(4)
and~(5) force $t=1$, $n_1=2$, as in~(6) and~(7) of the conclusion.
If $C$ lies in $\GO_n^+(q^2).2$ then again we must have $s=0$ and $C=L_d(2;0)$
as in~(8). Note that the conditions in (6), (7) and~(8) are mutually exclusive.

Now assume $r>0$. If $C$ lies in $\GO_{2ae}^\pm(q)\times\GO_{2r}^\pm(q)$
then we must have $s=r$ and thus $n_1>1$ since $C>\bC_G(x)$, as in~(1) or~(2)
of our conclusion. If $q=2$, $r=1$ and $C$ lies in $\Sp_{2n-2}(2)$ then in fact
we see $C$ lies in a subgroup
$\GO_{2n-2}^+(2)=(\GO_{2n-2}^+(2)\GO_2^+(2))\cap G$ (see
Proposition~\ref{prop:over Sp}(7)) and we are back in (1) or~(2). Finally, in
all of the cases discussed above apart from the
last one, $\Sub_G(x)$ must equal the relevant maximal subgroup, as can be seen
by using our earlier descriptions of maximal subgroups of classical groups
containing the normaliser of a Sylow $d$-torus in
Propositions~\ref{prop:over SL}, \ref{prop:over SU}, \ref{prop:over Sp},
\ref{prop:over SO}, \ref{prop:over SO+} and~\ref{prop:over SO-}.
\end{proof}

\begin{rem}   \label{rem:triality}
 Observe that for $n=4$, cases (5) and (8) of Theorem~\ref{thm:subn SOeven+}
 with $d=4$ are conjugate under triality (if $q$ is even). (See
 \cite[Tab.~8.50]{BHR}.) Since a Sylow $d$-torus $\bS_d$
 of $\bG=\SO_8$ can be chosen invariant under a triality automorphism $\tau$
 commuting with $F$, this shows that subnormalisers can behave strangely with
 respect to upward extensions, even by cyclic groups: if $x\in\bS_d^F$ has
 centraliser $L_4(2;0)$, for example, then its subnormaliser in $G=\SO_8^+(q)$
 is as given in (6), (7) or~(8) of Theorem~\ref{thm:subn SOeven+}, but its
 subnormaliser in $\hat G:=G\langle\tau\rangle$ is $\hat G$. A similar
 phenomenon for 3-elements in $\tw2F_4(2)'<\tw2F_4(2)$ was already observed in
 \cite[\S5.3]{Ma25}.
\end{rem}

We complete our investigations by considering subnormalisers in orthogonal
groups of minus-type.

\begin{thm}   \label{thm:subn SOeven-}
 Let $G=\SO_{2n}^-(q)$ with $n\ge4$, let $\ell\nmid q$ be a prime such that
 Sylow $\ell$-subgroups of $G$ are abelian and $\bS_d\le\bG$ a Sylow $d$-torus
 where $d=e_\ell(q)$. Set $e:=e_\ell(q^2)$ and let $a,r$ be as defined above.
 Then for $x\in\bS_d^F$ an $\ell$-element we have one of
 \begin{enumerate}[\rm(1)]
  \item $\Sub_G(x)=\bN_G(\bS_d)$ if $\bC_G(x)=\bC_G(\bS_d)$;
  \item $\Sub_G(x)=\big(\GO_{2ae}^+(q)\times\GO_{2r}^-(q)\big)\cap G$ if
   $2ae/d$ is even and $\bC_{G}(x)=L_d(n_1,\ldots,n_t;r)$ with $n_1>1$;
  \item $\Sub_G(x)=\big(\GO_{2ae}^-(q)\times\GO_{2r}^+(q)\big)\cap G$ if
   $2ae/d$ is odd and $\bC_{G}(x)=L_d(n_1,\ldots,n_t;r)$ with $n_1>1$;
  \item $M=(\GO_{2e}^-(q)\wr\fS_a)\cap G$ if $r=0$, $a$ is odd, $d$ is even
   and $\bC_{G}(x)=L_d(1,\ldots,1;e)$; or
  \item $\Sub_G(x)=G$ otherwise.
 \end{enumerate}
\end{thm}

\begin{proof}
The argument is similar to but easier than the one for the orthogonal groups
of plus-type. With the usual reductions we may assume $a\ge2$ and $C:=\bC_G(x)$
is a $d$-split Levi subgroup of type $L_d(n_1,\ldots,n_t;s)$ lying in one of
the maximal subgroups $M$ in (1)--(4) of Proposition~\ref{prop:over SO-}.
If $r=0$, so we are in Case~(4), the containment $C\le M$ forces $s=e$ and
$n_1=\ldots=n_t=1$, so we reach~(4) of the conclusion. If $r>0$ then $C\le M$
gives $s=r$ and we are in (2) or (3) of the conclusion; note that again
Case~(3) of Proposition~\ref{prop:over SO-} does not appear as a subnormaliser
since in that case $C$ is contained in a proper subgroup
$\GO_{2n-2}^-(2)=(\GO_{2n-2}^-(2)\GO_2^+(2))\cap G$ of $\Sp_{2n-2}(2)$ (see
Proposition~\ref{prop:over Sp}(8)).
\end{proof}

\subsection{An algebraic group construction}   \label{subsec:alg}

Looking back on the results obtained for the various types of classical groups
we observe that, as in the case of exceptional groups, the subnormalisers of
semisimple $\ell$-elements (in abelian Sylow subgroups) occurring in classical
groups are always normalisers of suitable subsystem subgroups, very similar to
the situation for semisimple elements in algebraic groups in
\cite[Thm~6.8]{Ma25}. We do not see, though, how to deduce this \emph{a priori}.
Nevertheless, there is a construction in the ambient algebraic group closely
related to subnormalisers in the finite case.

Let $\bG$ be connected reductive in characteristic~$p$ and $F:\bG\to\bG$ a
Frobenius map with respect to an $\FF_q$-structure.
Let $d\ge1$ and $\bT_d\le\bG$ a Sylow $d$-torus with relative Weyl group
$W_d:=\bN_G(\bT_d)/\bC_G(\bT_d)$, where $G:=\bG^F$.

\begin{prop}   \label{prop:algsubnorm}
 Let $\bL\le\bG$ be a $d$-split Levi subgroup of $\bG$ containing $\bT_d$. Then
 $$\bS:=\langle\bL^w\mid w\in W_d\rangle$$
 is an $F$-stable connected reductive subgroup of $\bG$ and $\bL$ is a
 $d$-split Levi subgroup of $\bS$. Furthermore, $\bS$ is the connected
 component of 
 $$\hbS:=\hbS(\bL,\bT_d):=\langle\bN_\bG(\bT_d)^F,\bL\rangle$$
 and $\bN_G(\bS)=\hbS^F$.
\end{prop}

\begin{proof}
Observe that $\bZ^\circ(\bL)_{\Phi_d}$ is a central $d$-torus of~$\bL$, hence
contained in $\bT_d\le\bL$ and so
$\bC_\bG(\bT_d)\le\bC_\bG(\bZ^\circ(\bL)_{\Phi_d})=\bL$. Thus, $\bL^w$ is
well-defined for $w\in W_d$. It is clear from their definition
that both $\hbS,\bS$ are closed and $F$-stable, and $\bS$ is connected, being
generated by conjugates of the connected group $\bL$ (see
\cite[Prop.~1.16]{MT11}). Furthermore, as $\bN_G(\bT_d)=\bC_G(\bT_d).W_d$
and $\bC_\bG(\bT_d)\le\bL$, $\bS$ is normalised by $\hbS$, and the factor
group is a quotient of $W_d$, hence finite. It follows by
\cite[Prop.~1.13]{MT11} that indeed $\bS$ is the connected component of the
identity of $\hbS$. Since $\bL$ is $d$-split in $\bG$, we have
$\bL=\bC_\bG(\bZ^\circ(\bL)_{\Phi_d})$, and so
$\bL=\bC_{\bS}(\bZ^\circ(\bL)_{\Phi_d})$ is also $d$-split in $\bS$.

To show that $\bS$ is reductive, let $\bT\le\bL$ be a maximal torus containing
$\bT_d$. Observe that by
\cite[Prop.~26.4]{MT11} we may consider $W_d$ as a subquotient of the Weyl
group of~$\bT$; we denote by $\hat W_d$ the full preimage. Let $\Phi$ be the
root system of~$\bG$ with respect to $\bT$ and let $\bU_\al\le\bG$,
$\al\in\Phi$, denote the corresponding root subgroups. 
Let $\Psi$ be the smallest $p$-closed subsystem of $\Phi$ containing
$\{\al^w\mid \bU_\al\le\bL, w\in \hat W_d\}$ (see \cite[Def.~13.2]{MT11}). Then
$\bH:=\langle\bT,\bU_\al\mid\al\in\Psi\rangle$ is connected reductive by
\cite[Thm~13.6]{MT11}. Since
$$\bL=\langle\bT,\bU_\al\mid\bU_\al\le\bL\rangle$$
by \cite[Thm~8.17(g)]{MT11}, we clearly have $\bS\le\bH$. On the other
hand, if $\bU_\al,\bU_\beta\le\bS$ and $\al+\beta\in\Psi$ then also
$\bU_{\al+\beta}\le\bS$ by the commutator formulas, and since
$\bU_\al\le\bL$ if and only if $\bU_{-\al}\le\bL$, inductively we also get
$\bU_{-\al-\beta}\le\bS$. Hence, for all $\bU_\al\le\bS$ we
have $\langle\bU_\al,\bU_{-\al}\rangle\le\bS$. Since this contains
representatives for the reflection $s_\al$ corresponding to $\al$, we also find
$\bU_{s_\al(\beta)}\le\bS$ for all $\bU_\beta\le\bS$. So by the
definition of $p$-closed we have $\bU_\al\le\bS$ for all $\al\in\Psi$.
Hence $\bS=\bH$, and so $\bS$ and $\hbS$ are reductive.

Finally, if $g\in\bN_G(\bS)$ then $\bT_d^g$ is a Sylow $d$-torus of $\bS$. So
there is $h\in\bS^F$ with $gh\in\bN_G(\bT_d)$, and thus
$g\in \bS^F\bN_G(\bT_d)$, showing $\bN_G(\bS)\le\hbS^F$, whence $\bN_G(\bS)=\hbS^F$.
\end{proof}

Now assume $\bG$ is simple and $\ell\ne p$ is a prime such that Sylow
$\ell$-subgroups of $G:=\bG^F$ are abelian, and $d=e_\ell(q)$. Then $\bT_d$
contains a Sylow $\ell$-subgroup of~$G$ by \cite[Prop.~2.4]{Ma14}, and all
centralisers of $\ell$-elements of $\bG^F$ are $d$-split Levi subgroups
of~$\bG$.

So for $x\in\bT_d^F$ an $\ell$-element, $\bL:=\bC_\bG(x)$ is a $d$-split Levi
subgroup of $\bG$ containing $\bT_d$ and we can construct the associated
group $\hbS(x):=\hbS(\bL,\bT_d)$ as in Proposition~\ref{prop:algsubnorm}. Note
that if $\bT_d'$ is another Sylow $d$-torus of $\bG$ containing $x$ then
$\bT_d,\bT_d'\le\bL$ are Sylow $d$-tori of $\bL$, hence $\bL^F$-conjugate, and
so are $\bN_G(\bT_d),\bN_G(\bT_d')$. Thus, $\hbS(x)$ is independent of the
choice of $\bT_d$ and only depends on $x$.

\begin{cor}
 In the above setting we have $\Sub_G(x)\le \hat S(x):=\hbS(x)^F$ and
 $\hbS(x)\le\Sub_\bG(x)$.
\end{cor}

\begin{proof}
Let $P\le G$ be a Sylow $\ell$-subgroup of $G$ contained in $\bT_d^F$, so
containing $x$. As $P$ is abelian by assumption, we have
$\Sub_G(x)=\langle\bN_G(P),\bC_G(x)\rangle$ by \cite[Prop.~2.12]{Ma25}.
Since $\bN_\bG(\hbS)^F$ is equal $\bN_G(P)$ by \cite[Prop.~2.4]{Ma14} this
shows the first inclusion. The second follows from the description of
$\Sub_\bG(x)$ in \cite[Thm~6.8]{Ma25}.
\end{proof}

\section{Subnormalisers in symmetric groups}   \label{sec:symm}
In this section we describe subnormalisers of $p$-elements in symmetric
groups $\fS_n$ with abelian Sylow $p$-subgroups. The results are very
similar to those for the special linear groups but easier to show. We write
$n=ap+r$ with $0\le r<p$, where have $a\le p-1$ as otherwise the Sylow
$p$-subgroups of $\fS_n$ are non-abelian. Note that a Sylow $p$-subgroup
of~$\fS_n$ is then elementary abelian, and the conjugacy class of a
$p$-element in $\fS_n$ is uniquely determined by its number of cycles of
length~$p$.

\begin{prop}   \label{prop:symm}
 Let $G=\fS_n$, $p$ a prime and $n=ap+r$ with $a,r\le p-1$. Let
 $x\in G$ have cycle type $(p)^k$. Then
 $$\Sub_G(x)=\begin{cases}
   \big(C_p.C_{p-1}\big)\wr\fS_a\times \fS_r& \text{if $k=a$ ($x$ is picky)},\\
   \fS_p\wr\fS_a& \text{if $r=0$ and $k=a-1\ge1$},\\
   G& \text{otherwise}.
 \end{cases}$$
 If $p>2$ and so $x\in\fA_n$ then $\Sub_{\fA_n}(x)=\Sub_G(x)\cap\fA_n$.
\end{prop}

\begin{proof}
A Sylow $p$-subgroup $P$ of $\fS_n$ has normaliser
$\big(C_p.C_{p-1}\big)\wr\fS_a\times \fS_r$, while the centraliser of
$x$ has the form $C_p\wr\fS_k\times\fS_{n-kp}$. Thus, if $x\in P$
and $k=a$ then $\bC_G(x)\le\bN_G(P)$ and hence $x$ is picky by
\cite[Cor.~2.7 and Prop.~2.12]{Ma25}. Now assume $k<a$. If $r=0$ and $k=a-1$
then clearly, $\bC_G(x)$ and $\bN_G(P)$ are both contained in a subgroup
$M=\fS_p\wr\fS_a$, and in fact, $M$ is generated by these, so we are in
case~(2) of the conclusion. So finally assume $k<a-1$, or $r>0$ and $k=a-1$.
Then $\bC_G(x)$ contains a symmetric group $\fS_{2p}$, respectively
$\fS_{r+p}$, and $\Sub_G(x)=\langle\bC_G(x),\bN_G(P)\rangle$ acts
transitively on $\{1,\ldots,n\}$, which forces $\Sub_G(x)=G$ by Jordan's
theorem.

Now assume $p>2$. Clearly, $\Sub_{\fA_n}(x)\le\Sub_G(x)\cap\fA_n$. For the
converse, we may assume $p>3$ by explicit computation, and then the same
reasoning as for $\fS_n$ applies.
\end{proof}

Note that for $n=6$, $p=3$, the first two cases in Proposition~\ref{prop:symm}
are conjugate under the exceptional outer automorphism.

The fact that $x$ is picky for $k=a$ was already shown by Mar\'oti, Mart\'inez
Madrid and Moret\'o \cite[Thm~3.9]{MMM}. The situation for primes
$p\le\sqrt{n}$ is much more involved and many different types of subnormalisers
can arise. For the prime~2 they were completely determined by Mart\'inez Madrid
\cite{MM25}. For example in $\fS_{15}$ there are eight different classes of
subnormalisers of 2-elements.

\section{On subnormalisers in sporadic groups}   \label{sec:spor}

The Sylow $p$-subgroups in sporadic simple groups $G$ are cyclic of prime
order in most cases, and then the subnormaliser of any non-trivial
$p$-element $x\in G$ is just $\bN_G(\langle x\rangle)$, by
\cite[Prop.~2.12]{Ma25}.
Here we consider the remaining instances of abelian Sylow subgroups:

\begin{prop}   \label{prop:spor}
 Let $G$ be a sporadic simple group and $p$ a prime such that Sylow
 $p$-subgroups of $G$ are abelian, but not cyclic of prime order.
 Let $x\in G$ be a $p$-element. Then $\Sub_G(x)=G$ unless $(G,p)$ are
 as given in Table~\ref{tab:spor}.
\end{prop}

\begin{table}[htb]
\caption{Subnormalisers in sporadic groups}   \label{tab:spor}
$\begin{array}{c|ccc}
 G& p& |\bC_G(x)|& \Sub_G(x)\\
\hline
  M_{11}& 3&  18& \text{picky}\\
     J_2& 5&  50& \text{picky}\\
     Suz& 5& 300& J_2.2\\
 Fi_{22}& 5& 600& \Om_8^+(2).\fS_3\\
\end{array}$
\end{table}

\begin{proof}
The sporadic groups $G$ satisfying the assumptions of the proposition are
$J_1$ for $p=2$, $M_{11},M_{22},M_{23},HS$ for $p=3$,
$J_2,Suz,He, Fi_{22},Fi_{23},Fi_{24}'$ for $p=5$, $Co_1,Th,B$ for $p=7$,
and the monster group for $p=11$. For the smaller groups the
claim is easily verified by a computer calculation. We just comment on the
larger cases (variations of the given arguments would in fact also allow to
treat the smaller cases by hand).
For $G=Fi_{24}'$ with $p=5$ there is a unique class of 5-elements. Now a
Sylow 5-subgroup of $G$ is contained in a maximal subgroup $M=Fi_{23}$, where
the subnormaliser of any 5-element $x$ equals $M$. Since
$|\bC_G(x)|>|\bC_M(x)|$ this shows $\Sub_G(x)=G$ by
Proposition~\ref{prop:Burnside}. For $G=Th$ with $p=7$
there is again a unique class of 7-elements. Here a Sylow 7-subgroup $P$ of
$G$ is contained in a maximal subgroup $M=\tw3D_4(2).3$, and this group has
one class of 7-elements $x$ with $\Sub_M(x)=M$ (see
Theorem~\ref{thm:subn exc}). Since $|\bN_G(P)|>|\bN_M(x)|$ we conclude
$\Sub_G(x)=G$ (using again Proposition~\ref{prop:Burnside}). For $G=B$ with
$p=7$ we use our previous result on the maximal subgroup $M=Th$ by noting
that centralisers of 7-elements in $G$ are larger than in $M$. Finally assume
$G$ is the monster group and $p=11$. Here $G$ has a unique class of
11-elements, but the centraliser of such an element has order not dividing the
order of a Sylow 11-normaliser, which is maximal in $G$, so again
$\Sub_G(x)=G$.   \par
Since $G=Fi_{22}$ is a maximal subgroup of $\tw2E_6(2)$, the subnormalisers of
the unique class of 5-elements $x\in G$ can be derived from
Theorem~\ref{thm:subn exc} for an upper bound and
Remark~\ref{rem:triality} for a lower bound.
\end{proof}

The subnormaliser $J_2.2$ in $Suz$ of course contains the subnormaliser~$J_2$
in the maximal subgroup $G_2(4)$ of $Suz$ found in Theorem~\ref{thm:subn exc}.

\begin{proof}[Proof of Theorem~\ref{thm:main}]
Let $S$ be quasi-simple. By \cite[Lemma~2.15]{Ma25} we may assume $\bZ(S)=1$,
so $S$ is simple. For $S$ a sporadic group, the subnormalisers were found in
Proposition~\ref{prop:spor}, while for $S=\fA_n$ they are known by
Proposition~\ref{prop:symm}, respectively \cite{MM25} when $p=2$.

So finally assume $S$ is simple of Lie type. If $p$ is the defining
characteristic of $S$, then $S\cong\PSL_2(p^f)$ and the non-trivial
$p$-elements are picky by \cite[Thm~6.1]{MMM} or \cite[Prop.~3.6]{Ma25}.
If $p$ is not the defining characteristic, then either $p>2$ or again
$S=\PSL_2(q)$. In the former case, subnormalisers are determined in
Sections~\ref{sec:exc} and~\ref{sec:class}. So let $S=\PSL_2(q)$ with
$q\equiv3,5\pmod8$. The cases when a 2-element $x\in S$ is picky are described
in \cite[Thm~6.1]{MMM} and \cite[Lemma~3.7]{MS25}. Otherwise, $\Sub_S(x)$
properly contains the normalisers $\fA_4$ of a Sylow 2-subgroup and of a maximal
torus of order $(q-1)/2$ or $(q+1)/2$ and hence equals $S$.
\end{proof}


\end{document}